\documentclass[12pt,a4paper]{article}

\usepackage[utf8]{inputenc}
\usepackage[english]{babel}

\usepackage[a4paper,margin=1in]{geometry}
\usepackage{tikz}
\usetikzlibrary{arrows.meta,positioning,fit,calc,backgrounds}
\usepackage{amsmath,amssymb,amsthm,mathtools,bm}
\usepackage{siunitx}
\usepackage{booktabs,tabularx,longtable,multirow,array}

\usepackage{graphicx}
\usepackage{float}
\usepackage{placeins}
\usepackage{adjustbox}
\usepackage{newunicodechar}
\newunicodechar{₂}{$_2$}
\usepackage{natbib}
\usepackage{pdflscape}
\usepackage[ruled,vlined]{algorithm2e}

\usepackage[font=small,labelfont=bf,justification=centering]{caption}
\usepackage{subcaption}

\usepackage[dvipsnames]{xcolor} % for OliveGreen etc.
\usepackage{soul}
\sethlcolor{lightgray}

\usepackage{tikz}
\usetikzlibrary{positioning,calc}
\usepackage{pgfplots}
\pgfplotsset{compat=1.18}
\usepackage{pgfplotstable}
\tikzstyle{block} = [draw,rectangle,thick,minimum height=2em,minimum width=2em,fill=green!10!white]

\usepackage{tikz}
\usetikzlibrary{arrows.meta,positioning,fit,calc}
\usepackage{enumitem}

\usepackage{authblk}
\usepackage{pifont} % for \ding

\usepackage[normalem]{ulem}
\usepackage{url}
\usepackage{setspace}
\usepackage{hyperref}
\hypersetup{colorlinks=true,linkcolor=blue,citecolor=blue,urlcolor=blue}

\usepackage{attachfile}

\title{\Large Integrating Hydrogen into Ontario’s Energy Hub: A Robust, Carbon-Aware Framework for Power–Heat–Transport} %\\ \hamed{Robust Carbon-Aware Co-Optimization of a Multi-Energy Hub with Hydrogen Under Carbon Tax and Net-Zero Policy Pathways: The Case of the Province of Ontario, Canada}}

\author[a]{Hossein Mirzaee}
\author[b]{Mostafa Mostafavi Sani\thanks{Corresponding Author}}
\author[c]{Hamed Samarghandi}
\affil[a]{\small{University of British Columbia, Vancouver, British Columbia, Canada, V6T 1Z4}}
\affil[b]{\small{Department of Industrial Engineering, Dalhousie University, Halifax, NS, Canada, B3H 4R2}}
\affil[c]{\small{Edwards School of Business, University of Saskatchewan, Saskatoon, SK, Canada, S7N 5A7}}

\providecommand{\keywords}[1]{\small\textbf{\textit{Keywords---}} #1}

\date{\vspace{-1.5cm}}

\begin{document}

\maketitle

\begin{abstract}
Decarbonizing electricity generation, heating, and transportation simultaneously requires integrated planning tools that can coordinate multiple energy production sources and demand points while remaining reliable under uncertainty. This paper develops a carbon-aware and uncertainty-resilient optimization framework for a grid-connected multiple source hub that co-optimizes electricity, heating, cooling, and transport energy services with an explicit hydrogen sub-hub. The proposed model is formulated as a mixed-integer linear program (MILP) over a 25-year planning horizon (2025--2050). The hub integrates renewable electricity (Photovoltaic and wind), dispatchable resources (including natural-gas-based conversion), storage systems, demand response, and a hydrogen subsystem comprising an electrolyzer and hydrogen storage to supply hydrogen-vehicle demand and provide temporal flexibility. Two policy archetypes are examined: a Carbon Tax (price instrument) implemented via an emissions cost term in the objective, and a Net-Zero pathway (quantity instrument) implemented through an emissions-trajectory constraint that tightens toward net-zero by 2050. To hedge feasibility-critical operational uncertainty, the deterministic model is extended using budgeted robust optimization (Bertsimas--Sim) with $\pm 30\%$ deviation bounds and a tunable uncertainty budget. The developed scheme is applied to the province of Ontario, Canada; the results indicate substantial long-term hydrogen expansion, with the electrolyzer capacity increasing from 300~MW (2025) to 3{,}800~MW (2050) and hydrogen storage from 2{,}000~MWh to 37{,}000~MWh (2050), accompanied by sharply higher hydrogen production. Compared with deterministic solutions, robust solutions preserve feasibility under adverse realizations but incur a moderate robustness premium of approximately 6.6--9.0\% in total cost across the policy cases studied, while slightly reducing hydrogen utilization and renewable share and increasing reliance on dispatchable balancing. Overall, the framework clarifies the reliability--cost--carbon trade-off faced by planners and highlights the value of additional low-carbon flexibility to mitigate the emissions impact of robustness.
\end{abstract}


\FloatBarrier

\section{Introduction}

The global energy system is undergoing a structural transition driven by climate mitigation imperatives and the need to reduce dependence on carbon-intensive fuels \citep{staffell2019role}. Fossil fuels have long underpinned economic growth, but are also the main source of energy-related greenhouse gas emissions, and still account for a large share of global electricity generation \citep{bogdanov2021low}. As a response, many jurisdictions are targeting net-zero scenarios and are increasing investments in cleaner sources of power, such as electricity generated from renewable sources, as well as new sources of energy production such as hydrogen. Renewable energy is expected to form a significant part of the global energy mix by the middle of the century, and hydrogen is increasingly seen as a potentially important carrier that can be generated from low-carbon power, stored for extended periods, and used as zero-carbon fuel for a range of applications \citep{khakimov2024hydrogen}.

In this context, integrated energy hubs have become a major concept for the integration of multiple-carrier energy systems \citep{mohammadi2018optimal}. An energy hub provides a unified platform for converting, storing, and distributing multiple forms of energy (e.g., electricity, heating and cooling, and fuels), enabling system-wide optimization rather than managing each carrier in isolation. By exploiting synergies between carriers, for example, using waste heat, converting surplus electricity into storable chemical energy, or shifting loads across vectors, energy hubs can improve efficiency, reduce losses, and facilitate greater penetration of renewables \citep{mohammadi2018optimal}. Recent reviews have noted that energy hubs are becoming a fundamental building block for sustainable energy systems to meet economic and environmental objectives \citep{papadimitriou2023comprehensive,Siroos2025}.

Transportation electrification and new alternative vehicle technologies are key aspects of this change and directly relate to energy hub development. Gasoline and diesel vehicles are predominant, especially in terms of overall numbers, but are also a major source of CO\textsubscript{2} emissions and local pollutants. Electric vehicles (EV), especially when powered by a clean energy source, are a lower-carbon alternative and, with growing numbers, represent a large and variable electrical load \citep{yu2022electric}. Electric vehicles also play a role in providing storage through vehicle-to-grid (V2G) technologies, providing additional flexibility to accommodate intermittent renewable energy sources. Additionally, hydrogen fuel cell vehicles (FCEVs) are also gaining popularity for transportation applications requiring long driving ranges and fast refueling, directly linking transportation electrification and decarbonization to hydrogen production and refueling infrastructure development \citep{staffell2019role}. Therefore, energy hubs increasingly need to consider transportation, electricity, and hydrogen within a single integrated framework.

However, the planning and operation of these types of multi-carrier hubs are further complicated due to the presence of pervasive uncertainties, which include the variability of the output from the renewable generators and the varying nature of the electricity, heat, and EV/FCEV demand, among others. In addition, long-term drivers, including the costs of technologies, demand, and climate policies such as carbon pricing, are subject to change over time. Hence, the failure to address these uncertainties can lead to underperforming or unreliable solutions, particularly when high penetration of renewables and flexible, mobile demand are involved in the solutions. Consequently, recent work on energy hubs and multi-energy systems underscores the need for optimization approaches that explicitly account for uncertainty while maintaining economic viability and environmental performance \citep{lasemi2022comprehensive}. This motivation is reinforced by evidence that uncertainty in costs and policy conditions can materially alter decarbonization pathways, particularly in integrated energy systems where carbon pricing affects technology choice and incentives for renewable deployment \citep{dilek2026carbon}.

Among the various energy carriers and storage possibilities that can be incorporated in an energy hub, hydrogen has been identified as having significant promise in terms of enhancing resilience and enabling deep decarbonization. The production of hydrogen from surplus renewable energy through electrolysis, long-term storage, and subsequent use in electricity production through fuel cells or direct use as a clean energy source in heating and transportation has been recognized as a promising approach \citep{khakimov2024hydrogen,staffell2019role}. The ability of hydrogen to offer long-term storage, including seasonal storage, has been recognized, allowing it to be utilized in conjunction with short-term storage options. Hydrogen has been identified to have the ability to absorb energy during times of high renewable energy production, thereby providing energy during prolonged low-renewable periods or peak demand times, allowing for greater penetration of renewables. At the same time, the use of hydrogen instead of fossil fuels has been recognized to offer environmental advantages in terms of reduced direct CO\textsubscript{2} and pollutant emissions \citep{staffell2019role}. Additionally, widespread use of hydrogen has been recognized to potentially alter electricity-gas interactions and price relationships, since hydrogen production can be made from both electricity and methane \citep{densing2026will}.

Robust optimization is particularly relevant in this setting because it seeks decisions that remain feasible under bounded adverse realizations of uncertain parameters, without requiring a fully specified probability distribution for every source of variability. This is especially attractive for integrated energy hubs, where renewable availability, multi-carrier loads, and EV/FCEV charging or refueling patterns can fluctuate considerably and may not be described reliably by a single probabilistic model. In addition, budgeted robust formulations allow the decision-maker to tune the level of conservatism, thereby balancing reliability and cost. Recent studies, therefore, show a growing use of robust and hybrid uncertainty-aware methods in energy-hub planning, and Canadian evidence suggests that explicitly modeling uncertainty can materially change both the optimal technology mix and the cost--emissions trade-off \citep{Siroos2025}.

Within this context, studies that \emph{jointly} treat hydrogen as an integral part of a multi-carrier energy hub, explicitly model both EV and FCEV fleets as interacting with the hub, and internalize greenhouse gas emissions through a carbon-cost mechanism under uncertainty are quite rare. Most studies focus on economic optimality without fully capturing environmental costs, or treat hydrogen, transportation, and uncertainty in a more simplified or decoupled manner. This paper contributes to filling this gap by developing a robust optimization model for a large-scale integrated energy hub that: (i) embeds a hydrogen sub-hub (electrolyzers, hydrogen storage, and fuel cells) alongside electricity, heating, and cooling subsystems; (ii) explicitly represents EV and FCEV loads and their interaction with the hub; (iii) internalizes emissions through an emissions-based cost term to reflect carbon pricing; and (iv) quantifies the robustness premium required to preserve feasibility under uncertain renewable output and demand conditions. Applied to the province of Ontario, Canada, the results show that robustness improves operational security but entails a moderate cost premium of approximately 6.6--9.0\% in the policy cases studied, while modestly reducing hydrogen utilization and renewable share and increasing reliance on dispatchable balancing resources. These findings clarify the reliability--cost--carbon trade-off faced by planners of integrated provincial energy systems.

The remainder of this paper is organized as follows. Section \ref{sec:LR} reviews the relevant literature and highlights the research gaps that motivate this work. Section \ref{sec:PD} describes the problem under study, including the structure of the proposed energy hub, key assumptions, and major modeling considerations. Section \ref{sec:MF} presents the mathematical formulation, including the objective function and the constraints of the model. Section \ref{sec:Results and Discussions} introduces the case study, describes the main data inputs, and presents the numerical results along with a detailed discussion of the insights obtained. Finally, Section \ref{sec:CC} concludes the article and outlines possible directions for future work.

\section{Literature Review}\label{sec:LR}

The concept of integrated energy hubs (EH) has been proposed as a major paradigm in the integration, storage, and distribution of electricity, heating, and fuels, which has been shown to potentially increase efficiency and reduce emissions relative to the management of each carrier in isolation \citep{Mohammadi2018,Papadimitriou2023}. At the same time, many deep decarbonization scenarios have recognized that, in addition to electrification, other carriers will be required if high penetration of renewables is to be realized cost-effectively \citep{staffell2019role}. One such carrier, which has re-emerged as a major enabler in this context, is hydrogen. Hydrogen has been recognized as a major long-term energy storage medium and carrier that complements electrification with many uses, including transportation, and therefore has the potential to reduce the consumption of fossil fuels in a structural sense \citep{staffell2019role}. All of these considerations will be relevant in the context of the Ontario grid, which will be examined in this study.

Among the various options that can be integrated into the EH, hydrogen is particularly interesting in providing long-duration balance, including seasonal balance, which has been supported through empirical evidence and modeling exercises, suggesting that hydrogen can absorb surplus energy when it is abundant, while providing balance during extended periods of low renewable penetration or winter peaks, which is not economically viable through battery technology \citep{Guerra2020,Petkov2020,Gabrielli2020}. At the same time, the use of hydrogen in generators, boilers, or vehicles offers clear benefits in reducing greenhouse gasses by eliminating direct pollutant emissions \citep{staffell2019role}. As a consequence, hydrogen emerges as a natural complement to the electric-dominated EHs, which justifies the explicit consideration of hydrogen within the planning model.

Electrification of transportation, including the use of alternative technologies, further supports the need to include integrated EHs in the formulation of the optimal solution, since electric vehicles are being deployed in large numbers, which in turn introduce new, flexible and sizeable electric loads, which can be utilized in turn as flexible storage units through vehicle-to-grid (V2G) technologies \citep{yu2022electric}. At the same time, hydrogen fuel cell electric vehicles are becoming increasingly popular in the transport sector, particularly in situations where range extension and rapid refueling are required, which are particularly critical in heavy-duty transport, further strengthening the interrelationship between decarbonization of transport and the provision of hydrogen infrastructure \citep{bhuiyan2025hydrogen,staffell2019role}. In practice, the need to support the transition from diesel-based electric vehicles to cleaner options, including EVs and FCEVs, while satisfying the residual diesel-based energy requirements provides the context within which the center under investigation operates, including the need to support the transition to cleaner options, including EVs and FCEVs, while satisfying the residual diesel-based requirements \citep{Mohammadi2017,Siroos2025}. This multi-modal, multi-carrier setting is the backdrop for the hub considered in this paper.

\subsection{Deterministic planning for multi-energy hubs}

Deterministic multi-period investment analysis studies have also provided insight into the potential for integrating renewables, distributed storage, and flexible conversion technologies that are available through EH and hybrid renewable energy systems (HRES). Representative models optimize multiple energy carriers and distributed storage over a planning horizon to match fluctuating supply and demand, given a point forecast \citep{Pizzuti2024, Tsao2023, Ling2024}. The structural benefits of such a multi-period optimization approach to investment planning are clear, although often a perfect foresight assumption is made, and costs are often assumed to be constant, limiting adaptability to falling costs and changing demand.

Hydrogen-related solutions, such as the concept ``Power to Gas'' and hydrogen storage, are being increasingly included in the deterministic planning literature. For example, some papers such as \cite{Alizad2022} and \cite{Sheng2024} examine the potential of using electrolyzers and hydrogen storage together with renewable energy sources to absorb generated energy and enable sector coupling. However, applications at the community or regional scale remain comparatively narrow: many works focus on industrial hubs or large integrated systems, leaving fewer insights for smaller, mixed-load systems where transportation, buildings, and local industry must be co-planned \citep{pizzuti2024novel,tsao2023distributed}. Similarly, hydrogen is often treated as an auxiliary \emph{add-on} technology that absorbs excess electricity rather than as a co-equal carrier with dedicated production, storage, and end-use services such as mobility, limiting the understanding of how long-duration storage and hydrogen refueling interact with electric loads in transport-coupled hubs \citep{wu2024collaborative}.

A related strand of deterministic work examines the interaction between EH and electrified transportation. Models incorporating EV charging and, in some cases, V2G highlight the potential of coordinated charging schedules to reduce costs and enable higher penetration of renewables \citep{yu2022electric}. However, these formulations typically either omit FCEVs or represent hydrogen demand in aggregate form, rather than explicitly modeling hydrogen production, storage, and refueling within the hub \citep{bhuiyan2025hydrogen,mohammadi2018optimal}. Consequently, they provide only a partial view of the infrastructure and operational requirements of dual-fleet (EV/FCEV) systems.

From a policy perspective, deterministic EH/HRES models often prioritize cost minimization and, at most, treat emissions as secondary constraints or separate objectives. Explicit carbon-pricing mechanisms that monetize emissions within the investment objective remain relatively rare, despite their importance in aligning planning decisions with climate policy targets \citep{Wang2021}. As a result, capacity expansions and technology choices in many deterministic studies may not fully reflect the incentives and trade-offs induced by carbon pricing or tightening emissions constraints \citep{petkov2020power, gabrielli2020seasonal}.

In summary, deterministic planning studies have established foundational insights into the benefits of multi-carrier coordination and the potential of hydrogen and EV integration, but they tend to (i) be tuned to large-scale systems with perfect foresight, (ii) treat hydrogen as an auxiliary carrier rather than a central, co-optimized subsystem, (iii) provide limited treatment of dual EV/FCEV fleets coexisting with diesel vehicles, and (iv) only partially internalize carbon policy instruments such as explicit carbon prices. These limitations motivate the move toward uncertainty-aware formulations.

\subsection{Dynamic planning under uncertainty}

To cope with variability and uncertainty, the literature has moved toward dynamic planning frameworks that explicitly model uncertain inputs and their temporal evolution. A common division of studies assigns multi-stage stochastic programming (MSSP) to long-term drivers such as demand growth, technology costs, and climate policy, as well as robust or distributionally robust optimization (RO/DRO) to short-term volatility in renewables and loads. MSSP captures sequential information revelation and enables adaptive investment recourse as scenarios unfold \citep{Ioannou2019,Hou2021,Liu2018}, while RO delivers operational reliability under worst-case intra-horizon realizations and has been applied in two-stage and multi-period settings, including in power and hydrogen-involved networks \citep{Tian2020,Zhang2023,Sun2023}.

Recent contributions have extended these approaches to more realistic representations of uncertainty in EHs and multi-energy systems. Scenario-based MSSP models explore how long-term uncertainty in technology costs, fuel prices, and policy instruments, such as carbon pricing, affect optimal capacity expansion and operational strategies. In parallel, RO and DRO formulations seek solutions that remain feasible and near-optimal in worst-case or ambiguity-aware distributions of uncertain parameters, thus prioritizing reliability in the face of volatile renewable output and demand \citep{Lasemi2022}. These methods are particularly relevant when integrating large shares of renewables and flexible loads such as EVs and FCEVs, whose availability and charging/refueling patterns are inherently uncertain and can significantly impact hub operations. Recent planning evidence further indicates that uncertainty in the transport-sector can bias both infrastructure sizing and operating cost when alternative-fuel vehicle adoption and refueling/charging requirements are modeled jointly with energy-system decisions, reinforcing the value of uncertainty-aware planning when EV and hydrogen mobility loads interact with supply systems \citep{najafzad2026two}.

Planning studies also examine staged transitions toward high-renewable systems, highlighting the importance of multi-stage policies and recourse actions. Multi-stage stochastic or robust frameworks allow investments to be updated as new information becomes available, improving the robustness of long-term decarbonization pathways \citep{Nunes2018,Santos2016}. However, even within these advanced formulations, hydrogen is frequently modeled in a simplified way, often as a generic storage option, without fully capturing its multi-sectoral role or its interaction with transportation demands. Recent works also emphasize the value of hedging and contracting strategies for green hydrogen producers under uncertainty, suggesting that uncertainty-aware planning is important not only for system operation, but also for hydrogen market participation and risk management \citep{palmer2025hedging}.

Despite these advances, three limitations persist. \textbf{(i) Coupling depth:} few models endogenize a \emph{hydrogen sub-hub} that co-optimizes electrolyzers, H$_2$ storage, and H$_2$ end-uses (e.g., fuel cells and refueling) \emph{together} with electricity/heat networks and detailed transport demands, including both EVs and FCEVs and their coexistence with legacy diesel fleets. \textbf{(ii) Objective completeness:} many MSSP/RO/DRO models optimize cost (and occasionally treat emissions as constraints or secondary objectives) but less frequently monetize carbon explicitly within the objective via carbon pricing, thus limiting sensitivity to policy signals and under-representing the economic value of emissions reductions \citep{Wang2021}. \textbf{(iii) Methodological balance:} pure RO can be conservatively biased, leading to potential over-sizing of assets, whereas pure MSSP can be scenario-heavy and computationally demanding; hybrid stochastic--robust constructs that exploit the strengths of both approaches remain under-explored in the context of multi-carrier hubs with hydrogen integration and dual-mode transportation. Complementary evidence from policy-focused studies indicates that incentives accelerating storage innovation and cost reduction can improve system flexibility, which can reduce the dependency on fossil-backed balancing in uncertainty-aware operations \citep{jin2026sustainably}.

\subsection{Limitations and contributions}

The present study addresses these limitations by:
\begin{enumerate}
    \item Embedding a hydrogen sub-hub that includes electrolyzers, H$_2$ storage, and fuel cells into a multi-carrier EH with transportation loads (EVs/FCEVs and residual diesel vehicles).

    \item Adopting a robustness-forward formulation to ensure reliable operation under short-term uncertainty in renewables and loads.

    \item Internalization of emissions through an explicit carbon-cost term in the objective function to reflect the full impact of carbon pricing.
\end{enumerate}
 
This design directly targets the reliability--cost--carbon trade space that previous deterministic and single-paradigm uncertain models have only partially captured, especially for community- and provincial-scale planning with deep hydrogen integration and transport coupling \citep{Papadimitriou2023,Lasemi2022,staffell2019role}.

\section{Problem Description}\label{sec:PD}

This paper considers a grid-connected multi-energy hub located in Ontario, Canada. The hub is supplied by the provincial electricity system and is coupled with local heating, cooling and hydrogen demands, as well as transport demands of electric vehicles (EVs) and hydrogen vehicles (HVs). The planning horizon is defined over $Y=\{1,\dots,25\}$ years. For each year, the operation of the system is represented by a 48-hour horizon $T=\{1,\dots,48\}$, in which a typical winter day ($t\leq 24$) and a typical summer day ($t>24$) are studied.

Upstream electricity is provided by disaggregated purchase variables associated with nuclear, hydro, gas, and biofuel units, each constrained by a fraction of its provincial capacity and by minimum generation shares. The total import $P^{\mathrm{buy}}_{y,t}$ and export $P^{\mathrm{sell}}_{y,t}$ of the grid are limited and are not allowed to occur simultaneously. Effective purchase costs are obtained from normalized generation costs, while emissions are evaluated by technology-specific
emission factors.

Within the hub, local conversion and storage technologies are installed. A gas-fired combined heat and power (CHP) unit and a gas boiler are modeled, with fuel consumption linked to electrical and thermal outputs via efficiency and lower heating value. A thermal storage system with state-of-charge, charging and discharging is included, subject to capacity limits and charge/discharge exclusivity. The cooling demand is supplied by an electric chiller, an absorption chiller, and cold storage, with algebraic relations based on performance coefficients.

A hydrogen sub-hub is also considered consisting of an electrolyzer and a hydrogen storage tank with bounds and charge/discharge limits; its balance equations ensure that hydrogen-vehicle demand is fully met by storage discharge. Cost and performance improvements in renewable and hydrogen technologies are represented through exogenous learning multipliers applied to their operation and maintenance (O\&M) terms over the 25-year period. 

\begin{figure}[!t]
\centering
\resizebox{\linewidth}{!}{%
\begin{tikzpicture}[
    font=\footnotesize\sffamily,
    line cap=round,
    line join=round,
    >={Latex[length=2.2mm]},
    node distance=8mm and 10mm,
    module/.style={
        draw=black!70,
        line width=0.7pt,
        rounded corners=2.5pt,
        minimum width=25mm,
        minimum height=8mm,
        align=center,
        inner sep=3.5pt,
        fill=white
    },
    source/.style={module, fill=gray!8},
    converter/.style={module},
    bus/.style={
        module,
        minimum width=28mm,
        minimum height=9mm,
        fill=blue!6,
        font=\footnotesize\sffamily\bfseries
    },
    storage/.style={
        module,
        rounded corners=7pt,
        dash pattern=on 3pt off 2pt,
        fill=gray!4
    },
    demand/.style={module, fill=gray!12},
    group/.style={
        draw=black!45,
        line width=0.8pt,
        rounded corners=4pt,
        fill=black!1,
        inner sep=7pt
    },
    glabel/.style={
        font=\footnotesize\sffamily\bfseries,
        text=black!75,
        fill=white,
        inner sep=1.5pt
    },
    elec/.style={->, line width=1.05pt, draw=blue!70!black, shorten >=1pt, shorten <=1pt},
    heat/.style={->, line width=1.05pt, draw=orange!85!black, shorten >=1pt, shorten <=1pt},
    cool/.style={->, line width=1.05pt, draw=teal!70!black, shorten >=1pt, shorten <=1pt},
    htwo/.style={->, line width=1.05pt, draw=green!55!black, shorten >=1pt, shorten <=1pt}
]

% -------------------------
% Sources and main electricity layer
% -------------------------
\node[source] (grid) {Electricity\\Grid};
\node[source, below=7mm of grid] (pv)   {PV};
\node[source, below=7mm of pv]   (wind) {Wind};

\node[bus, right=24mm of pv] (ebus) {Electricity\\Bus};
\node[storage, above=16mm of ebus] (bat) {Battery\\Storage};
\node[demand, below=16mm of ebus]  (ed)  {Electric\\Demand / EV};

% -------------------------
% Heat subsystem
% -------------------------
\node[converter, right=25mm of ebus, yshift=34mm] (hp)  {Heat\\Pump};
\node[converter, above=7mm of hp]                  (aux) {Aux.\ Heater};
\node[storage,   right=15mm of hp]                 (tes) {Thermal\\Storage};
\node[demand,    right=15mm of tes]                (hd)  {Heat\\Demand};

% -------------------------
% Hydrogen subsystem
% -------------------------
\node[converter, right=25mm of ebus] (ely) {Electrolyzer};
\node[storage,   right=15mm of ely]  (h2s) {H$_2$\\Storage};
\node[demand,    right=15mm of h2s]  (h2d) {FCEV\\Demand};
\node[converter, below=7mm of h2s]   (fc)  {Fuel\\Cell};

% -------------------------
% Cooling subsystem
% -------------------------
\node[converter, right=25mm of ebus, yshift=-34mm] (ch)  {Electric\\Chiller};
\node[storage,   right=15mm of ch]                  (ces) {Cooling\\Storage};
\node[demand,    right=15mm of ces]                 (cd)  {Cooling\\Demand};

% -------------------------
% Subsystem frames (background)
% -------------------------
\begin{scope}[on background layer]
    \node[group, fit=(grid)(pv)(wind)]    (srcframe) {};
    \node[group, fit=(aux)(hp)(tes)(hd)]  (hframe)   {};
    \node[group, fit=(ely)(h2s)(h2d)(fc)] (hyframe)  {};
    \node[group, fit=(ch)(ces)(cd)]       (cframe)   {};
\end{scope}

\node[glabel, anchor=south west] at ($(srcframe.north west)+(2mm,0.8mm)$) {Sources \& Grid};
\node[glabel, anchor=south west] at ($(hframe.north west)+(2mm,0.8mm)$)   {Heat Subsystem};
\node[glabel, anchor=south west] at ($(hyframe.north west)+(2mm,0.8mm)$)  {Hydrogen Subsystem};
\node[glabel, anchor=north west] at ($(cframe.south west)+(2mm,-0.8mm)$)  {Cooling Subsystem};

% -------------------------
% Electricity branching point
% -------------------------
\coordinate (etrunk) at ($(ebus.east)+(9mm,0)$);
\draw[elec] (ebus.east) -- (etrunk);
\fill[blue!70!black] (etrunk) circle (0.8pt);

% Sources to electricity bus
\draw[elec] (grid.east) -- ++(7mm,0) |- ([yshift=4mm]ebus.west);
\draw[elec] (pv.east) -- (ebus.west);
\draw[elec] (wind.east) -- ++(7mm,0) |- ([yshift=-4mm]ebus.west);

% Battery charge/discharge (single two-sided arrow)
\draw[<->, line width=1.05pt, draw=blue!70!black, shorten >=1pt, shorten <=1pt]
    (bat.south) -- (ebus.north);

% Electric demand
\draw[elec] (ebus.south) -- (ed.north);

% Electricity to subsystems
\draw[elec] (etrunk) -- (ely.west);
\draw[elec] (etrunk) |- (hp.west);
\draw[elec] (etrunk) |- (aux.west);
\draw[elec] (etrunk) |- (ch.west);

% Heat pathway
\draw[heat] (hp.east) -- (tes.west);
\draw[heat] (tes.east) -- (hd.west);
\draw[heat] (aux.south) -- ++(0,-4mm) -| (tes.north);

% Hydrogen pathway
\draw[htwo] (ely.east) -- (h2s.west);
\draw[htwo] (h2s.east) -- (h2d.west);
\draw[htwo] (h2s.south) -- ++(0,-3mm) -- (fc.north);

% Fuel-cell electricity back to lower-quarter point of right edge of bus
\draw[elec] (fc.west) -| ($(ebus.south east)!0.25!(ebus.north east)$);

% Cooling pathway
\draw[cool] (ch.east)  -- (ces.west);
\draw[cool] (ces.east) -- (cd.west);

% -------------------------
% Boxed legend
% -------------------------
\coordinate (l1s) at ($(srcframe.south west |- cframe.south)+(5mm,-14mm)$);
\draw[elec] (l1s) -- ++(9mm,0) coordinate (l1e);
\node[anchor=west] (l1t) at ($(l1e)+(2mm,0)$) {Electricity};

\coordinate (l2s) at ($(l1t.east)+(7mm,0)$);
\draw[heat] (l2s) -- ++(9mm,0) coordinate (l2e);
\node[anchor=west] (l2t) at ($(l2e)+(2mm,0)$) {Heat};

\coordinate (l3s) at ($(l2t.east)+(7mm,0)$);
\draw[cool] (l3s) -- ++(9mm,0) coordinate (l3e);
\node[anchor=west] (l3t) at ($(l3e)+(2mm,0)$) {Cooling};

\coordinate (l4s) at ($(l3t.east)+(7mm,0)$);
\draw[htwo] (l4s) -- ++(9mm,0) coordinate (l4e);
\node[anchor=west] (l4t) at ($(l4e)+(2mm,0)$) {Hydrogen};

\begin{scope}[on background layer]
    \node[
        draw=black!45,
        line width=0.7pt,
        rounded corners=3pt,
        fill=white,
        inner sep=4pt,
        fit=(l1s)(l1e)(l1t)(l2s)(l2e)(l2t)(l3s)(l3e)(l3t)(l4s)(l4e)(l4t)
    ] {};
\end{scope}

\end{tikzpicture}%
}
\caption{Carrier-flow schematic of the proposed multi-energy hub. \\Solid rectangles denote sources and conversion technologies, dashed nodes denote storage units, shaded rectangles denote end-use demands, and colored arrows indicate electricity, heat, cooling, and hydrogen flows.}
\label{fig:energyhub_schematic}
\end{figure}
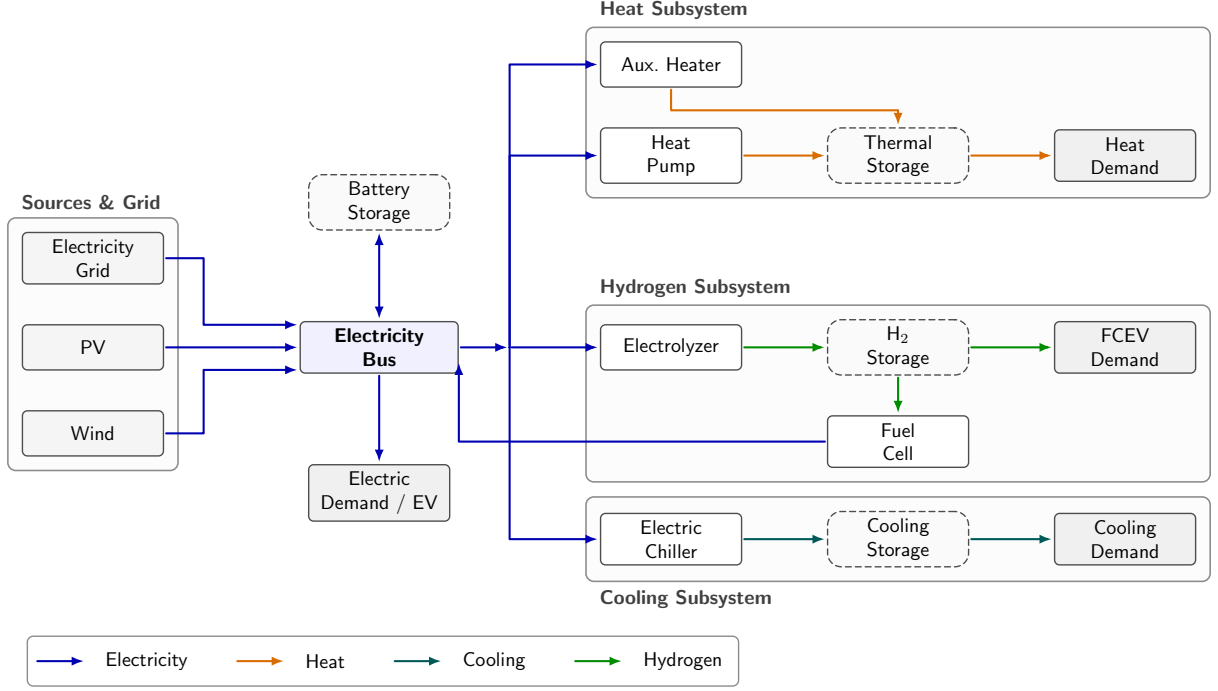

Figure \ref{fig:energyhub_schematic} illustrates the architecture and carrier interactions of the proposed multi-energy hub. The electricity bus acts as the central coupling point, receiving power from the grid and renewable sources (PV and wind), providing the electricity demand and the EV load, exchanging energy with battery storage, and feeding the main conversion technologies in the heating, cooling, and hydrogen subsystems. In the heating subsystem, the heat pump and the auxiliary heater provide thermal energy that can be stored in thermal storage and later dispatched to meet the heat demand. In parallel, the hydrogen subsystem converts electricity into hydrogen through the electrolyzer, stores it for later use, and supplies FCEV demand while also allowing reconversion to electricity through the fuel cell. The cooling subsystem is driven by electrical input and cooling storage to satisfy the cooling demand.

Demand response (DR) is incorporated for both electricity and heat. Upward and downward adjustments of the electric and thermal loads are allowed within prescribed fractions of the original demands, are controlled by binary variables, and are constrained to yield net-zero shifts over the 48-hour horizon. Penalties are associated with these adjustments and are included in the cost function.

To ensure analytical clarity and computational tractability, the proposed multi-energy hub model is developed under several assumptions. As described earlier, the proposed grid-connected energy hub in Ontario is modeled over a 25-year planning horizon, with annual operation represented by typical periods and aggregated time-dependent energy demands that are assumed to be fully satisfied in all operating periods. The technical characteristics of installed technologies are assumed to be known, whereas long-term cost variations are represented exogenously. The hydrogen subsystem is modeled through simplified conversion and storage relationships, with losses captured through efficiency parameters. EV and FCEV demands are represented in an aggregate form, without considering individual mobility behavior or charging congestion. The demand response is considered within predefined flexibility limits and is assumed to preserve net energy over the representative horizon. Finally, internal network constraints, forced outages, and short-term operational nonlinearities are neglected to maintain a tractable long-term planning model.

%========================================================
\subsection{Mathematical Formulation}\label{sec:MF}
%========================================================

A mixed-integer linear programming (MILP) framework is developed for the coordinated operation of the proposed multi-energy hub. The nominal version of this framework constitutes the deterministic benchmark, and the same structural model is later used as the basis for the robust counterpart. The operating horizon is represented by a set of years, and each year is discretized into a finite set of intra-year time intervals. Electricity, heat, cooling, and hydrogen carriers are jointly scheduled so that system demands are satisfied while operational limits, inter-temporal dynamics, and demand-response (DR) flexibility are respected. %The implementation details and the complete mathematical--computational mapping are documented in the supplementary implementation material. \hamed{\textbf{To Mostafa:} how can the readers access the supplementary implementation material?}

The deterministic backbone of the present formulation is adapted from the integrated energy-hub model in Ontario developed by \cite{Siroos2025}, which documents a realistic Ontario representation including PV, wind, nuclear, hydro, biofuel and natural-gas generation, EV charging/discharging interactions, and uncertainty treatment through CVaR-, IGDT- and robust-based analyses. In the present paper, the mentioned Ontario-specific operational backbone is retained as the nominal model, and is extended in three important directions:
\begin{enumerate}
    \item An explicit hydrogen sub-hub is added.
    \item Carbon-cost monetization is embedded as an objective.
    \item A robust optimization layer is introduced for uncertain parameters critical to feasibility.
\end{enumerate}
To avoid repetition, only the principal equations needed for the present analysis are summarized below, while the detailed baseline model can be found in \cite{Siroos2025}. Appendix \ref{Appendix:Notations} explains the notation used to formulate the problem.

\subsubsection{Objective functions}

A total cost objective is adopted so that economically efficient operating schedules are produced while accounting for energy trading, O\&M, emissions monetization, demand-response incentives/penalties, and storage cycling effects. The annual cost is aggregated across intra-year periods, and the total objective is obtained by summation over all years:
\begin{align}
\min \; J_{\text{cost}}
&:= \sum_{y\in\mathcal{Y}}\Bigg[
\sum_{t\in\mathcal{T}}
\Big(
\pi^{\text{buy}}_{y,t} P^{\text{buy}}_{y,t}
-\pi^{\text{sell}}_{y,t} P^{\text{sell}}_{y,t}
\Big) 
\;+\; C_{\text{cyc}}(y)
\;+\; C_{\text{OM}}(y)
\;+\; C_{\text{EM}}(y)
\;+\; C_{\text{DR}}(y)
\Bigg], \label{eq:obj_cost}
\end{align}
where $C_{\text{DR}}(y)=C_{\text{EDR}}(y)+C_{\text{HDR}}(y)$ is used to capture DR penalties/incentives in electricity and heat, and $C_{\text{EM}}(y)$ is used to monetize emissions through a carbon price. In this way, emission reductions are encouraged without introducing nonlinearity, and the resulting schedules remain comparable under different policy assumptions.

To explicitly represent environmental performance, a second objective is stated as the minimization of direct emissions from fuel use.
\begin{align}
\min \; J_{\text{emis}}
&:= \sum_{y\in\mathcal{Y}}\sum_{t\in\mathcal{T}}
P_f\Big(G^{\text{gas}}_{y,t}+G^{\text{bio}}_{y,t}\Big) . \label{eq:obj_emis}
\end{align}
This objective is included to enable transparent reporting of trade-offs between economic and environmental targets, particularly when carbon prices are uncertain or when an emissions policy is enforced as a hard constraint.
\begin{align}
\min \; J_{\text{cost}} + J_{\text{emis}}, \label{eq:weighted_sum}
\end{align}

\subsubsection{Constraints}

\paragraph{1) Electricity balance and carrier coupling.}
A nodal electricity balance is imposed in each period so that supply and demand are matched and the feasibility is preserved at the operating timescale. Electricity imports, renewable injections, and CHP electricity are required to cover electrical demands, conversion loads, and electricity exports:
\begin{align}
P^{\text{buy}}_{y,t} + P^{\text{PV}}_{y,t} + P^{\text{WT}}_{y,t} + P^{\text{ECHP}}_{y,t}
=
P^{\text{EL}}_{y,t} + P^{\text{EV}}_{y,t} + P^{\text{sell}}_{y,t}
+ P^{\text{ice}}_{y,t} + P^{\text{ec}}_{y,t} + P^{\text{ely}}_{y,t},
\quad \forall y,t.
\label{eq:el_balance}
\end{align}
The electrolyzer term is included to ensure that hydrogen production is endogenously linked to electricity availability and price conditions, thereby enabling power-to-hydrogen flexibility.

\paragraph{2) Heat and cooling balances.}
A heat balance is imposed so that heat production and thermal discharge collectively satisfy the heat demand and any thermal charging requirements, while also supplying thermally driven cooling when applicable:
\begin{align}
P^{\text{HCHP}}_{y,t} + P^{\text{HB}}_{y,t} + P^{\text{dch,H}}_{y,t}
=
L^{h}_{y,t} + P^{\text{ch,H}}_{y,t} + P^{\text{Hac}}_{y,t},
\quad \forall y,t.
\label{eq:heat_balance}
\end{align}
In this representation, the heat input from the absorption chiller is explicitly accounted for to preserve consistency between the coupled heat–cooling conversions. Cooling and auxiliary conversion relations (e.g. electric chiller, ice storage, absorption chiller) are imposed to ensure that cooling demand is met through physically consistent transformations. In this manner, inter-carrier substitutability (electric-driven versus heat-driven cooling) is captured within a unified framework.

\paragraph{3) CHP and boiler input--output relations.}
Fuel-to-energy conversion constraints are imposed to ensure that electric and thermal outputs remain consistent with fuel inputs and conversion efficiencies. Linear relations are used so that the resulting optimization remains MILP:
\begin{align}
P^{\text{HB}}_{y,t}   &= \frac{\text{LHV}\,\eta^{\text{B}}}{G^{\text{B}}_{y,t}}, \label{eq:boiler}\\
P^{\text{ECHP}}_{y,t} &= \frac{\text{LHV}\,\eta^{\text{CHP}}_e}{G^{\text{CHP}}_{y,t}}, \label{eq:chp_e}\\
P^{\text{HCHP}}_{y,t} &= \frac{\text{LHV}\,\eta^{\text{CHP}}_h}{G^{\text{CHP}}_{y,t}}. \label{eq:chp_h}
\end{align}
The upper bounds on fuel inputs are further imposed to reflect the equipment capacity limits, thus preventing infeasible operating points and ensuring that dispatch remains within the admissible region.

\paragraph{4) Heat storage dynamics and operating limits.}
Inter-temporal coupling is introduced through storage dynamics so that shifting of thermal energy across time can be optimally scheduled. For all $y$ and $t>1$, the heat storage state is updated as:
\begin{align}
H_{y,t}
=
H_{y,t-1}
+
\eta^{\text{ch}}_{\text{H}} P^{\text{ch,H}}_{y,t}
-
\frac{1}{\eta^{\text{dch}}_{\text{H}}} P^{\text{dch,H}}_{y,t},
\quad \forall y,\; t>1.
\label{eq:heat_storage_dyn}
\end{align}
An initial condition at $t=1$ is imposed to anchor the state trajectory within each year. Storage bounds are enforced to represent physical energy capacity limits:
\begin{align}
H^{\min}\le H_{y,t}\le H^{\max}, \quad \forall y,t. \label{eq:heat_state_bounds}
\end{align}
Charge/discharge power limits are enforced and are activated through binary mode selection:
\begin{align}
0\le P^{\text{ch,H}}_{y,t}\le P^{\max}_{\text{ch,H}} K^{\text{ch}}_{y,t},\;
0\le P^{\text{dch,H}}_{y,t}\le P^{\max}_{\text{dch,H}} K^{\text{dch}}_{y,t},
\quad \forall y,t, \label{eq:heat_power_bounds}\\
K^{\text{ch}}_{y,t}+K^{\text{dch}}_{y,t}\le 1,\quad
K^{\text{ch}}_{y,t},K^{\text{dch}}_{y,t}\in\{0,1\},
\quad \forall y,t.
\label{eq:heat_mutual}
\end{align}
Mutual exclusivity is imposed to prevent simultaneous charging and discharging, which would otherwise create artificial cycling and potentially distort cost signals.

In addition, a yearly net-zero condition is imposed so that the storage state does not drift across years when each year is modeled independently.
\begin{align}
\sum_{t\in\mathcal{T}} P^{\text{ch,H}}_{y,t} = \sum_{t\in\mathcal{T}} P^{\text{dch,H}}_{y,t},\quad \forall y.
\label{eq:heat_net_zero}
\end{align}
This condition is adopted to ensure comparability of annual schedules and to prevent the exploitation of end-effects when the year-end state is not explicitly linked across years.

\paragraph{5) Hydrogen subsystem: storage dynamics, coupling, and electrolyzer limits.}
Hydrogen production and storage are represented to capture long-duration flexibility and sector coupling. For all $y$ and $t>1$, the hydrogen inventory is updated as follows:
\begin{align}
H2_{y,t}
=
H2_{y,t-1}
+
\eta^{\text{ch}}_{\text{H2}} P^{\text{ch,H2}}_{y,t}
+
\eta^{\text{ely}} P^{\text{ely}}_{y,t}
-
\frac{1}{\eta^{\text{dch}}_{\text{H2}}} P^{\text{dch,H2}}_{y,t},
\quad \forall y,\; t>1.
\label{eq:h2_dyn}
\end{align}
An initial condition at $t=1$ is imposed for the same reason as in thermal storage, namely to define a well-posed trajectory. Inventory bounds are enforced as follows:
\begin{align}
H2^{\min}\le H2_{y,t}\le H2^{\max}, \quad \forall y,t. \label{eq:h2_state_bounds}
\end{align}
Charge/discharge limits and mutual exclusivity are enforced through binary mode selection to prevent simultaneous charging and discharging.
\begin{align}
0\le P^{\text{ch,H2}}_{y,t}\le P^{\max}_{\text{ch,H2}} K^{\text{ch,H2}}_{y,t},\;
0\le P^{\text{dch,H2}}_{y,t}\le P^{\max}_{\text{dch,H2}} K^{\text{dch,H2}}_{y,t},
\quad \forall y,t, \label{eq:h2_power_bounds}\\
K^{\text{ch,H2}}_{y,t}+K^{\text{dch,H2}}_{y,t}\le 1,\quad
K^{\text{ch,H2}}_{y,t},K^{\text{dch,H2}}_{y,t}\in\{0,1\},\quad \forall y,t.
\label{eq:h2_mutual}
\end{align}
Hydrogen use is coupled to hydrogen-serving demand so that the feasibility of the hydrogen service is guaranteed:
\begin{align}
P^{\text{H2,use}}_{y,t} = P^{\text{HV}}_{y,t}, \qquad
P^{\text{H2,use}}_{y,t} = P^{\text{dch,H2}}_{y,t},
\quad \forall y,t.
\label{eq:h2_meet}
\end{align}
The electrolyzer is constrained by its nameplate capacity:
\begin{align}
0\le P^{\text{ely}}_{y,t} \le P^{\max}_{\text{ely}},\quad \forall y,t.
\label{eq:ely_cap}
\end{align}
These constraints are introduced so that hydrogen production is coordinated with electricity availability and so that hydrogen delivery is ensured through the storage–discharge pathway, thus representing a consistent power-to-hydrogen-to-service chain.

\paragraph{6) Demand response bounds and logical constraints.}
The demand response is modeled as bounded upward and downward deviations around baseline demands. Limits are imposed as fractions of the baseline demand to reflect comfort or service-quality restrictions, and binary variables are introduced to prevent simultaneous upward and downward actions:
\begin{align}
0\le P^{\uparrow}_{el,y,t} \le MR^{\uparrow}_{el} L^{el}_{y,t} I^{\uparrow}_{el,y,t},\quad
0\le P^{\downarrow}_{el,y,t} \le MR^{\downarrow}_{el} L^{el}_{y,t} I^{\downarrow}_{el,y,t},
\quad \forall y,t, \label{eq:dr_el_bounds}\\
I^{\uparrow}_{el,y,t}+I^{\downarrow}_{el,y,t}\le 1,\quad
I^{\uparrow}_{el,y,t},I^{\downarrow}_{el,y,t}\in\{0,1\},
\quad \forall y,t. \label{eq:dr_el_logic}
\end{align}
Analogous constraints are imposed for heat DR. In this way, short-term flexibility is represented without allowing unrealistic oscillations, while economic incentives are captured through $C_{\text{DR}}(y)$ in the objective.

\paragraph{7) Ramping constraints.}
Inter-temporal ramping limits are imposed for selected generators so that physically plausible trajectories are produced and abrupt changes in output are prevented:
\begin{align}
-P^{\max}_{r} \le P_{y,t}-P_{y,t-1}\le P^{\max}_{r},\quad \forall y,\; t>1,
\label{eq:ramp_generic}
\end{align}
where $P^{\max}_{r}$ is set according to the corresponding unit (e.g. $ P_{\text{nuc}}$ or $ P_{\text{hyd}}$). Ramping constraints are included to ensure that dispatch solutions remain implementable and to avoid dispatch artifacts that may otherwise be produced by purely static balance constraints.

\paragraph{8) Domain constraints.}
All flow variables are constrained to be non-negative; all mode-selection variables are constrained to be binary. These domain restrictions are required to preserve physical meaning (e.g. negative imports or negative charging are disallowed) and to correctly enforce logical exclusivity.

\medskip
Equations \eqref{eq:obj_cost}--\eqref{eq:ramp_generic} collectively define a MILP in which multi-carrier balances, conversion physics, DR flexibility, and storage inter-temporal dynamics are coherently integrated. As a result, operating schedules are generated that are both physically feasible and economically interpretable under the adopted cost and emissions accounting framework.

The above equations define the common analytical core used throughout the numerical study. In the next section, this model is first solved with nominal inputs to establish the deterministic benchmark for Ontario and to identify the main operating and policy trends implied by the formulation. The same structure is then extended to the robust setting by protecting feasibility-critical uncertain parameters, making it possible to quantify the resulting robustness premium and to examine how uncertainty reshapes system cost, emissions, renewable utilization, and hydrogen deployment under Carbon Tax and Net-Zero policy settings.

\section{Results and Discussions}
\label{sec:Results and Discussions}

This section reports results for two implementations of the proposed Ontario energy-hub MILP: (i) a deterministic model solved under nominal (forecast) inputs and (ii) a robust optimization (RO) variant that accounts for uncertainty in selected operational parameters using a budgeted uncertainty set. The model operates on two coupled temporal layers. At the strategic level, investment and capacity decisions, such as electrolyzer sizing, hydrogen storage expansion, and technology deployment, are optimized over a 25-year planning horizon (2025--2050), capturing long-term decarbonization trajectories under carbon tax and net-zero policy pathways. At the operational level, each year's day-to-day system operation is represented by two representative 24-hour periods, a typical cold (winter) day and a typical warm (summer) day, yielding a 48-hour intra-year operating window. This two-day representation is chosen to capture the dominant seasonal variability in renewable availability, heating and cooling demands, and transport loads while keeping the model computationally tractable across the full planning horizon. The discussion proceeds in four steps: 
\begin{enumerate}
    \item Deterministic results are presented to establish baseline operational behavior and long-term investment trajectories.

    \item A solver-based sensitivity analysis identifies the parameters that most influence total cost and those that can threaten feasibility under large deviations.

    \item These feasibility-critical parameters motivate the uncertainty set in the robust formulation.

    \item Robust results are interpreted and compared with deterministic outcomes to quantify the robustness premium and derive managerial insights.
\end{enumerate}

Table~\ref{tab:parameters} summarizes the key input parameters introduced in the present model beyond the baseline established in \cite{Siroos2025}. The baseline Ontario hub parameters, including generation O\&M costs, technology efficiencies, thermal and cooling storage limits, demand response ratios, and EV-related parameters, are adopted directly from \cite{Siroos2025} and are not repeated here for brevity. Consequently, table \ref{tab:parameters} focuses on three sets of parameters that are new to the present model: the hydrogen sub-hub, the carbon policy instruments, and the robust optimization layer.

\begin{table}[H]
\centering
\caption{Key parameters introduced in the present model} 
\label{tab:parameters}
\begin{adjustbox}{width=\textwidth,center}
\begin{tabular}{llll}
\toprule
\textbf{Parameter} & \textbf{Symbol} & \textbf{Value} & 
\textbf{Source} \\
\midrule
\multicolumn{4}{l}{\textit{Hydrogen sub-hub}} \\
Electrolyzer capacity (2025) 
    & $P^{\max}_{\text{ely}}$ 
    & 300~MW 
    & \cite{ontario2022hydrogen} \\
Electrolyzer efficiency 
    & $\eta^{\text{ely}}$ 
    & 0.75 
    & \cite{Petkov2020} \\
H$_2$ storage capacity (2025) 
    & $H2^{\max}$ 
    & 2{,}000~MWh 
    & \cite{huang2025ontario} \\
H$_2$ storage minimum state 
    & $H2^{\min}$ 
    & 200~MWh (10\% of $H2^{\max}$) 
    & \cite{Guerra2020,Petkov2020} \\
H$_2$ storage charging efficiency 
    & $\eta^{\text{ch}}_{\text{H2}}$ 
    & 0.95 
    & \cite{Gabrielli2020} \\
H$_2$ storage discharging efficiency 
    & $\eta^{\text{dch}}_{\text{H2}}$ 
    & 0.95 
    & \cite{Gabrielli2020} \\
Max H$_2$ storage charge power 
    & $P^{\max}_{\text{ch,H2}}$ 
    & 200~MW 
    & \cite{Guerra2020} \\
Max H$_2$ storage discharge power 
    & $P^{\max}_{\text{dch,H2}}$ 
    & 200~MW 
    & \cite{Guerra2020} \\
\midrule
\multicolumn{4}{l}{\textit{Carbon policy}} \\
Carbon tax -- base year (2025) 
    & $\tau_{\text{CO}_2}$ 
    & \$80/tonne~CO$_2$e 
    & \cite{canada2025carbontax} \\
Carbon tax -- escalation (2025--2030) 
    & -- 
    & +\$15/tonne/year 
    & \cite{canada2025carbontax} \\
Carbon tax -- long-run (2030--2050) 
    & -- 
    & \$170/tonne~CO$_2$e (held constant) 
    & \cite{canada2025carbontax} \\
Net-zero emissions trajectory 
    & -- 
    & Linear reduction to zero by 2050 
    & Policy scenario \\
\midrule
\multicolumn{4}{l}{\textit{Robust optimization}} \\
Uncertainty deviation bound 
    & $\hat{a}_{ij}$ 
    & $\pm$30\% of nominal 
    & \cite{bertsimas2004price} \\
Uncertainty budget 
    & $\Gamma$ 
    & Varied (see Fig.~\ref{fig:cost_vs_gamma}) 
    & \cite{bertsimas2004price} \\
Protected constraints 
    & $\mathcal{I}^R$ 
    & Electricity and H$_2$ balance constraints 
    & This paper \\
\bottomrule
\end{tabular}
\end{adjustbox}
\end{table}

The deterministic and robust optimization models were coded in the General Algebraic Modeling System (GAMS) and solved with the CPLEX solver. CPLEX was used because the proposed formulation is an MILP and CPLEX is well-suited for large-scale MILP problems. All runs were performed on a personal laptop with a 12th Gen Intel\textsuperscript{\tiny\textregistered} Core\textsuperscript{\tiny TM} i5-1235U processor (1.30~GHz) and 16~GB of RAM (3200~MT/s). The results reported in the following section correspond to the optimal solution returned by CPLEX for the specified model instances.

\subsection{The deterministic model}
\label{subsec:deterministic}

\subsubsection{Operational decisions over the representative 48-hour horizon}
The deterministic operational solution is reported on a representative 48-hour horizon consisting of two 24-hour windows (interpreted as cold- and warm-day profiles). Table~\ref{tab:gen_energy_24h} summarizes hourly imports by generation segment (hydro, gas, bio), allocation of transport across the EV/HV/fossil pathways, electrolyzer electricity consumption, and the state of hydrogen storage.

Overall, total imports decrease from the first 24-hour window (typical cold day) to the second (typical warm day), along with a significan reduction in the share of natural gas compared to hydro. Transport allocation also increases for electrified and hydrogen-based pathways, where the EV share increases and the share of hydrogen-vehicles increases on a warm day compared to a cold day, along with a decrease in the share of fossil vehicles. Electrolyzer electricity consumption also increases on a warm day, and hydrogen storage has a net positive charge over both 24-hour periods, indicating that hydrogen is being accumulated rather than being fully utilized within a single 24-hour period. This is consistent with an expectation where the H$_2$ inventory accumulates within a given time period to meet future demand.

\begin{table}[!ht]
\centering
\caption{Energy supplied by source over the representative days}
\label{tab:gen_energy_24h}
\begin{tabular}{ccc}
\toprule
\textbf{Source} & \textbf{Cold day (MWh)} & \textbf{Warm day (MWh)} \\
\midrule
Hydro (hydropower) & 9{,}200  & 7{,}800  \\
Wind               & 4{,}800  & 5{,}200  \\
Solar PV           & 1{,}600  & 3{,}200  \\
Biofuel            & 1{,}200  & 1{,}200  \\
Natural Gas        & 6{,}800  & 3{,}200  \\
Hydrogen  & 2{,}626.4 & 3{,}637.2 \\
\midrule
\textbf{Total energy including hydrogen} & \textbf{26{,}226.4} & \textbf{24{,}237.2} \\
\bottomrule
\end{tabular}
\end{table}

\begin{figure}[!ht]
\centering
\includegraphics[width=\linewidth]{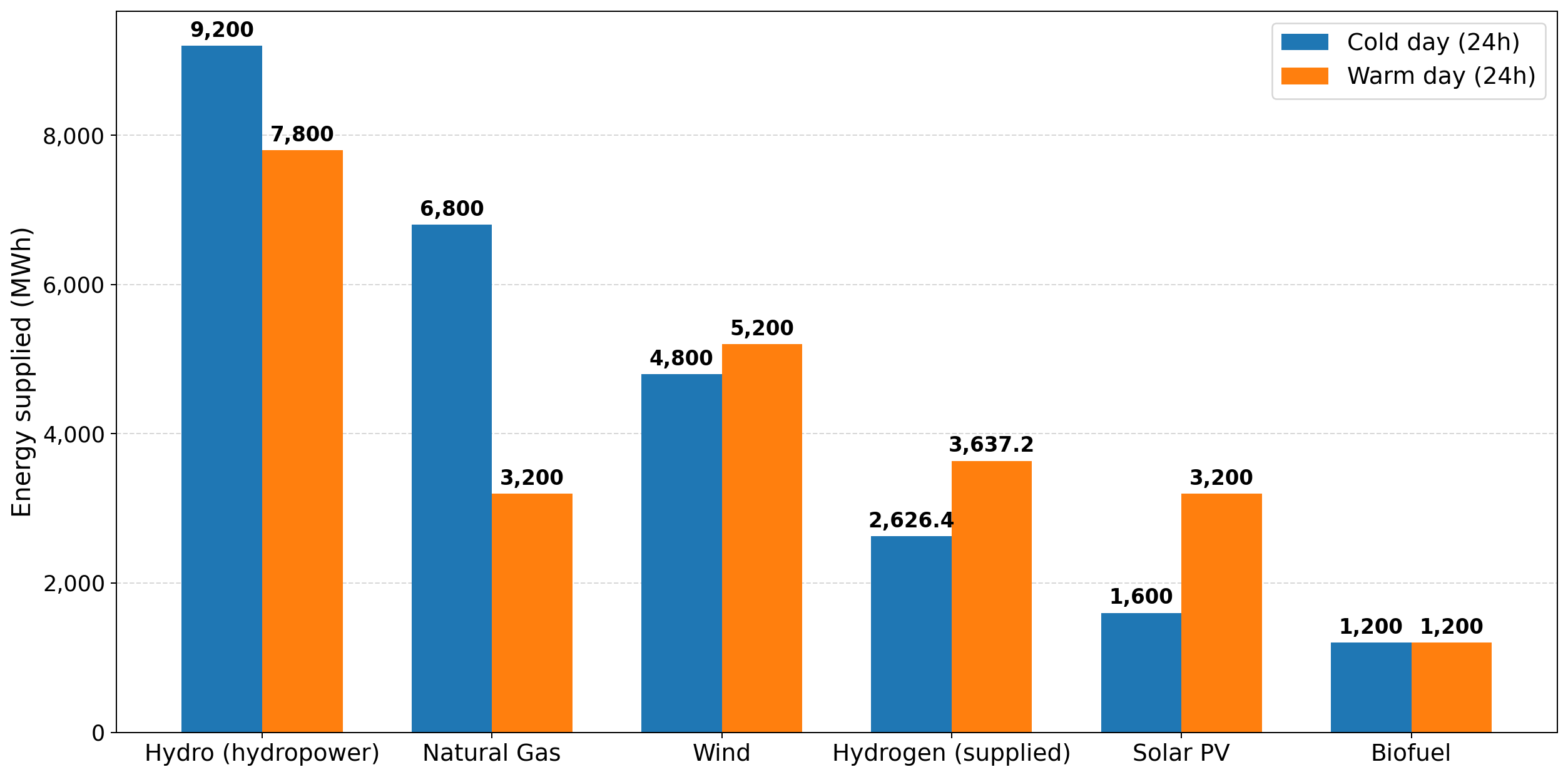}
\caption{Energy supply by source during the representative days}
\label{fig:energy_grouped_rep_days}
\end{figure}

\begin{figure}[H]
\centering
\includegraphics[width=\linewidth]{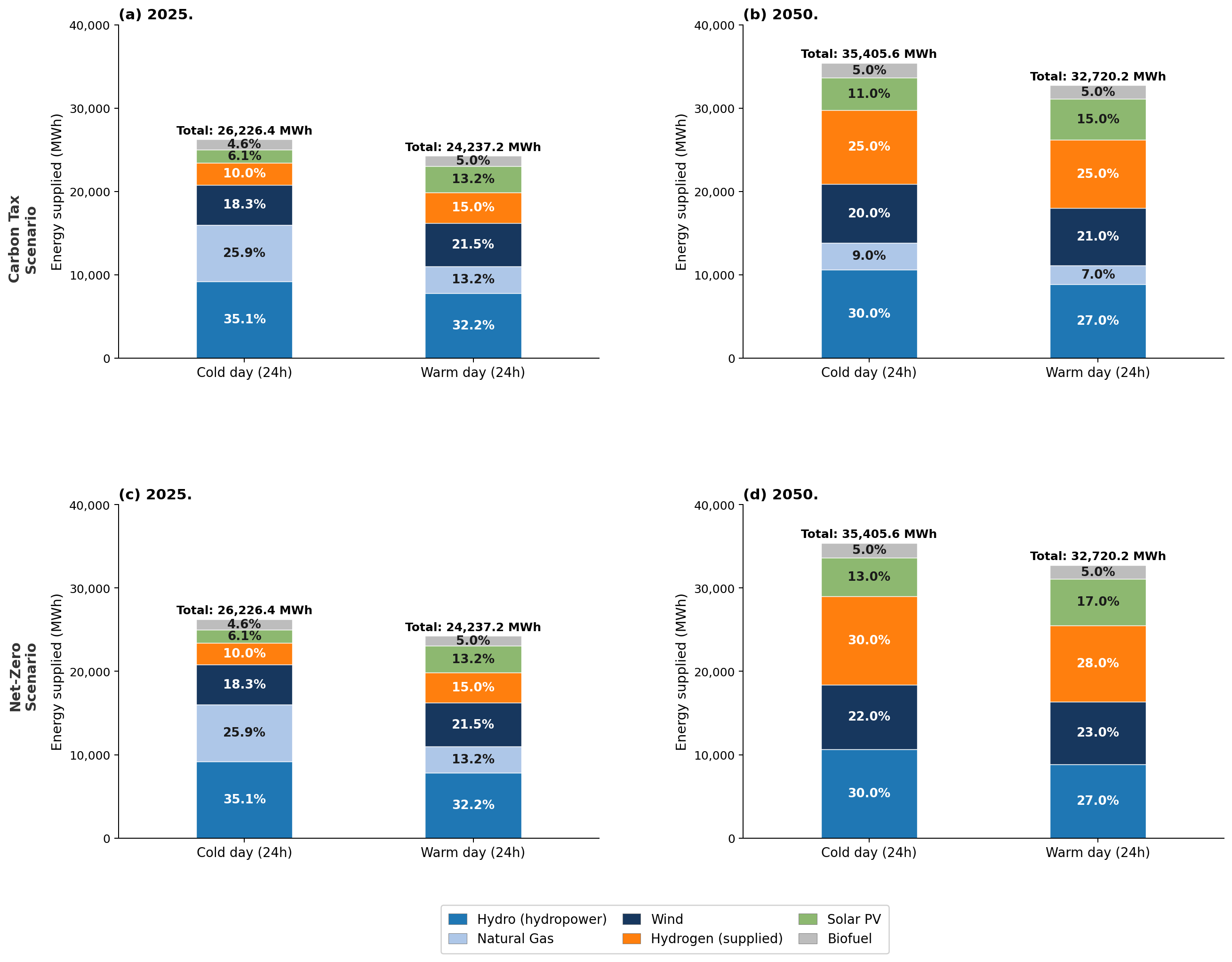}
\caption{Energy supply composition over representative cold and warm days at the beginning (2025) and end (2050) of the 25-year planning horizon, shown separately for the carbon tax scenario (panels a--b) and the net zero scenario (panels c--d).}
\label{fig:fig3_ab}
\end{figure}

The results in Table~\ref{tab:gen_energy_24h} and Figures \ref{fig:energy_grouped_rep_days} and \ref{fig:fig3_ab} show the typical intra-year operational snapshot used by the model: a 24-hour cold-day profile (for winter) and a 24-hour warm-day profile (for summer). These values represent the energy supplied over the two specific days within the 48-hour operating period. Although these values do not reflect the averages over the entire 2025--2050 planning horizon, Table~\ref{tab:gen_energy_24h} and Figures \ref{fig:energy_grouped_rep_days} and \ref{fig:fig3_ab} depict how the energy hub distributes the supply between carriers and technologies during different seasons.

There are several prominent patterns. First, more dispatchable resources, particularly natural gas, are required to support the cold day. On a warm day, the contribution of natural gas is 3,200 MWh; on a cold day, it is 6,800 MWh. This change corresponds to higher peak and thermal loads during winter, when dispatchable generation keeps the hub feasible and reliable. On a cold day, the contribution of hydro energy also increases from 7,800 MWh to 9,200 MWh. This suggests that when the net demand is higher, flexible, low-carbon supply plays a bigger role.

Furthermore, on warm days, the share of variable renewables increases. Wind increases from 4,800 to 5,200 MWh, while photovoltaic increases from 1,600 to 3,200 MWh. The seasonal availability of the representative profiles is reflected in these modifications. The fact that biofuel does not change between the two days indicates that it either has a fixed supply or acts as a reliable dispatch contributor.

In addition, hydrogen is present in both operating windows and acts as an energy supply source at the hub level (measured in MWh$_{H2}$). The hydrogen supply is higher on warm days (3,637.2 MWh$_{H2}$) than on cold days (2,626.4 MWh$_{H2})$, which is consistent with increased electrolyzer operation during improved electrical conditions, such as increased renewable output and decreased heating demand. This implies that hydrogen serves as a carrier of supportive energy. Through the hydrogen sub-hub, H$_{2}$ production and availability can fluctuate according to seasonal conditions, helping to meet transportation and other hydrogen-based demands.

\subsubsection{Long-term hydrogen build-out (2025--2050)}

\begin{figure}[ht]
\centering
\includegraphics[width=.8\textwidth]{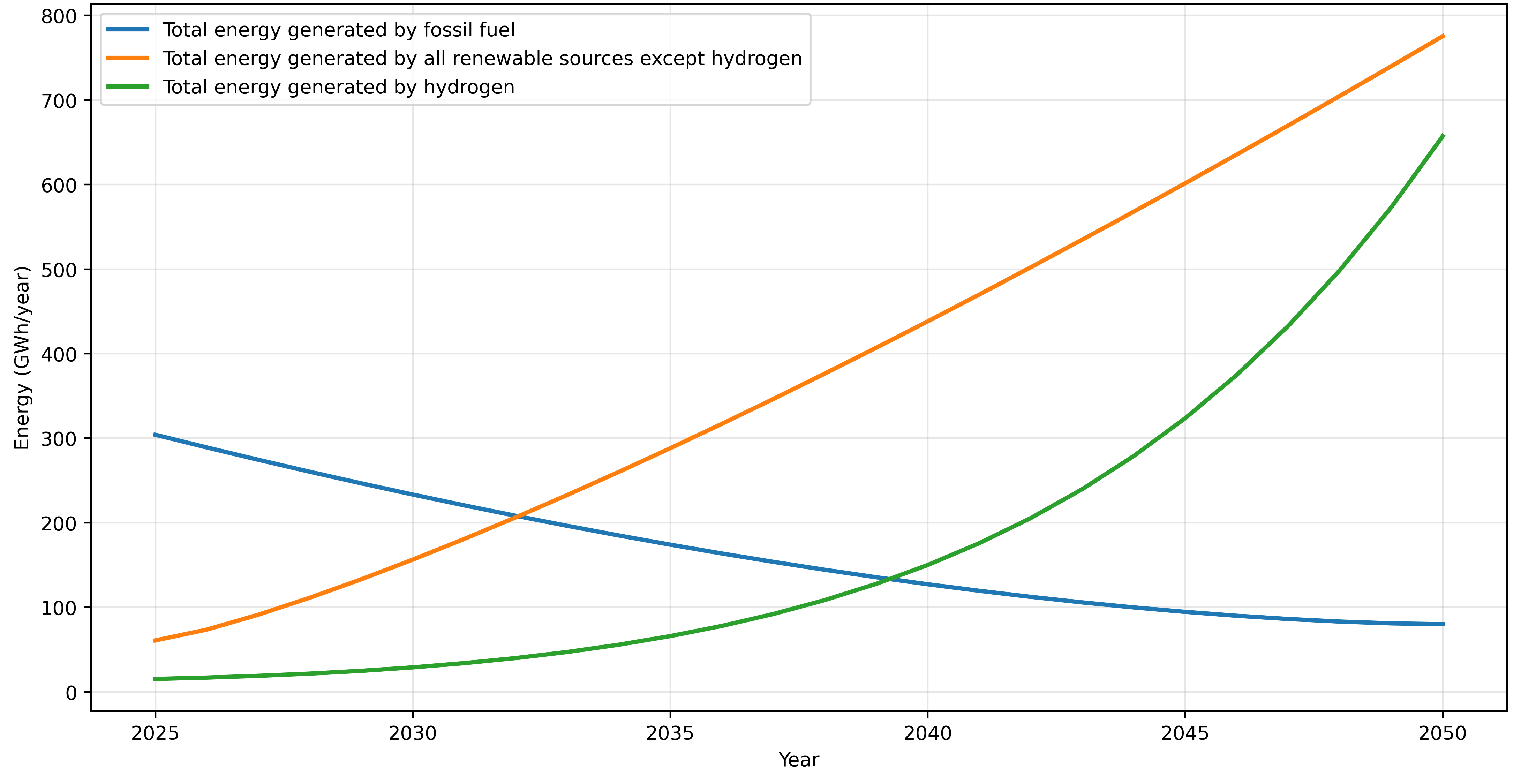}
\caption{Transition under the carbon tax pathway: fossil fuel energy use decrease while both electricity and hydrogen supplies increase over 2025--2050}
\label{fig:transition_elec_h2_fossil}
\end{figure}

Figure~\ref{fig:transition_elec_h2_fossil} shows the changes in the system as the carbon tax increases over the planning horizon. In 2025, fossil fuel energy accounts for approximately 304 GWh/year. This amount is roughly five times the electricity supplied, which is 60.8 GWh/year, and about twenty times the hydrogen supplied, which is 15.2 GWh/year. As carbon costs increase with time, the use of fossil fuels gradually decreases. It drops from 304 GWh/year in 2025 to approximately 233.3 GWh/year in 2030, 174.0 GWh/year in 2035, and 127.2 GWh/year in 2040. By 2050, it reaches a much lower level of about 80 GWh/year.

In parallel, both electricity and hydrogen supplies increase to meet the growing demand, but hydrogen grows more rapidly. Hydrogen supply increases from 15.2 GWh/year in 2025 to 28.9 GWh/year in 2030, 150.0 GWh/year in 2040, and 657.2 GWh/year in 2050. Electricity also increases steadily, going from 60.8 GWh/year in 2025 to 156.4 GWh/year in 2030, 438.2 GWh/year in 2040, and 775.5 GWh/year in 2050. Electricity remains the main energy carrier throughout this period because it is already widely used and can support more end uses with lower additional infrastructure costs. Hydrogen grows rapidly due to the targeted investment in electrolyzer conversion and hydrogen storage, which allows energy to be shifted over time and used in the future.

By 2050, the electricity supply is slightly above the hydrogen supply at 775.5 GWh/year compared to 657.2 GWh/year, making it about 18\% larger. This aligns with a system focused on electrification, supported by a rapidly expanding hydrogen sector, and a much smaller (but still present) role for fossil fuels under the carbon tax pathway.

\subsubsection{Deterministic policy effect: carbon tax vs net-zero}

Table \ref{tab:policy_det} outlines the difference in the result of two policy scenarios. The carbon tax scenario uses a price-based policy in which the optimization function is modified to include the cost of carbon. This allows optimization to balance the reduction of emissions with the cost. On the other hand, the Net-Zero scenario is based on a quantity/constraint-based policy. The optimization is constrained to reduce greenhouse gas emissions to net zero by 2050, which limits the solution space and restricts the development of technology.

\begin{figure}[H]
\centering
\includegraphics[width=\linewidth]{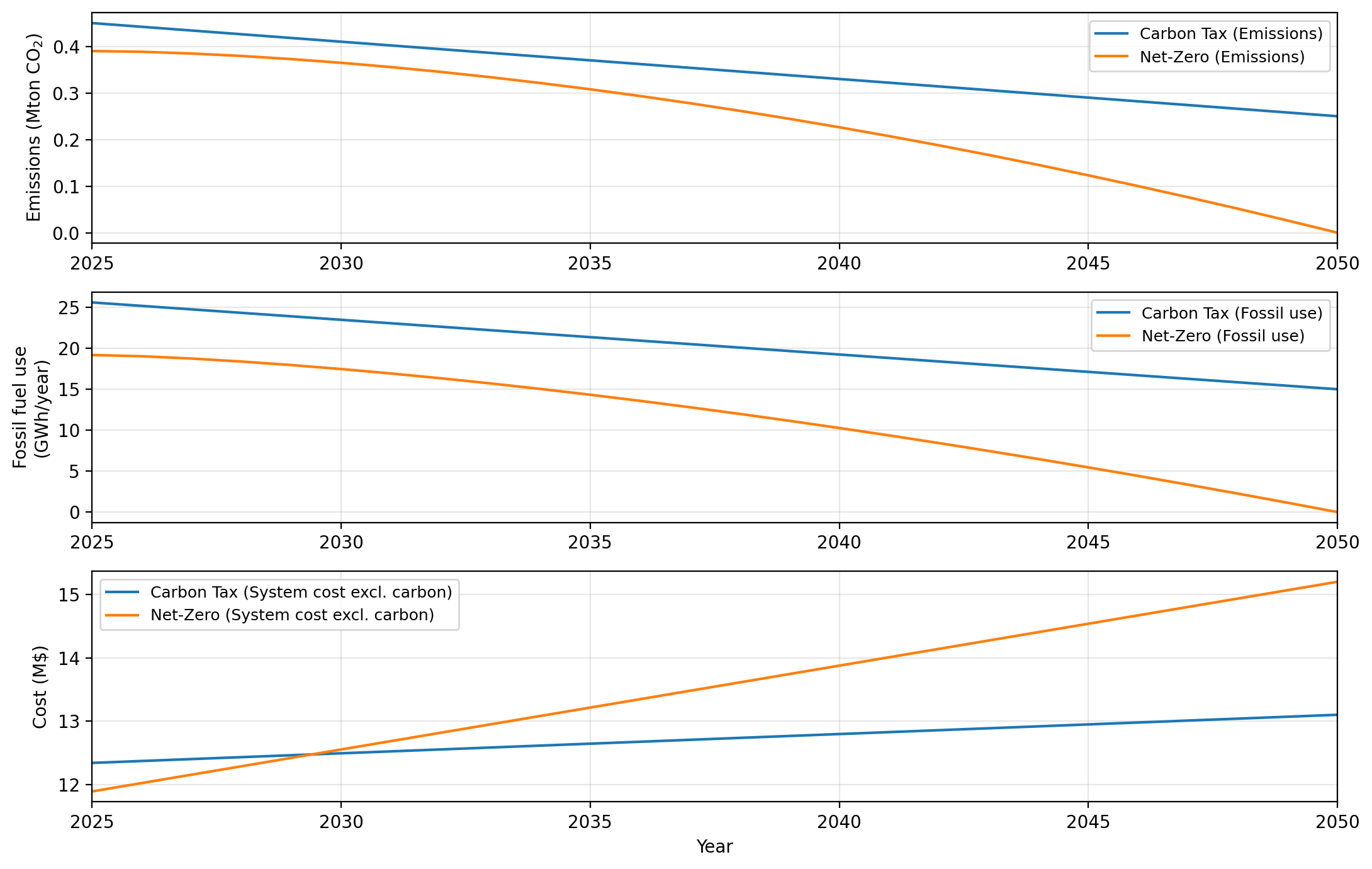}
\caption{Schematic long-run trajectories (2025--2050) illustrating the expected directional effects of carbon tax versus net zero policies: (top) emissions decline under both policies, with net zero enforcing a steeper reduction toward near-zero by 2050; (middle) fossil fuel use declines in both cases, with net zero requiring deeper phase-down by 2050; (bottom) system cost, excluding carbon costs, increases more under net zero due to tighter decarbonization requirements.}
\label{fig:policy_schematic}
\end{figure}

In the net zero scenario, greenhouse gas emissions are reduced to 0.39 Mton CO$_2$ compared to 0.45 Mton CO$_2$ in the carbon tax scenario. This is a reduction of 13.3\%. Similarly, in the net zero scenario, the consumption of fossil fuels is reduced to 19.18 GWh/year compared to 25.6 GWh/year in the carbon tax scenario. This is a reduction of 25.1\%. However, in the net zero case, hydrogen fuel production increases to 18.5 GWh/year compared to 15.2 GWh/year in the carbon tax scenario, which is an increase of 21.7\%. Furthermore, in the net zero case, the percentage of renewable fuels increases to 35.2\% compared to 28.5\% in the carbon tax scenario, which is an increase of 6.7\%. Figure \ref{fig:policy_schematic} is an indexed representation of the trends. As shown, the net zero policy reduces greenhouse gas emissions and fossil fuel consumption while increasing low-carbon fuel supplies.

Interestingly, the total cost also reduces from 12.34 to 11.89 million dollars (M\$), a 3.65\% decrease, under the net zero policy. This result fits an accounting structure where the carbon tax case leads to higher direct carbon payments, while the net zero case mainly achieves emissions reductions through technology changes and operational shifts. In simpler terms, with a sufficiently strong constraint, the system can reduce significant carbon costs by directly lowering emissions, which balances some of the necessary system adjustments. The main takeaway is that these two methods promote decarbonization differently: carbon tax encourages gradual substitution based on cost, while net zero necessitates deeper structural changes toward renewables and hydrogen to meet the emissions target.

\begin{table}[H]
\centering
\caption{Deterministic policy comparison: net zero vs. carbon tax. \\Positive values indicate net zero exceeds carbon tax. \\Negative values indicate reductions under net zero.}
\label{tab:policy_det}
\begin{adjustbox}{width=\textwidth,center}
\begin{tabular}{lcccc}
\toprule
\textbf{Metric (Daily)} & \textbf{Carbon Tax (Det)} & \textbf{Net-Zero (Det)} & \textbf{$\Delta$ (NZ--CT)} & \textbf{$\Delta$\%} \\
\midrule
Total Cost (M\$) & 12.34 & 11.89 & -0.45 & -3.65\% \\
Total Emissions (MTon CO$_2$) & 0.45 & 0.39 & -0.06 & -13.33\% \\
Fossil Fuel Use (GWh/year) & 25.60 & 19.18 & -6.42 & -25.08\% \\
H$_2$ Production (GWh/year) & 15.20 & 18.50 & +3.30 & +21.71\% \\
Renewable Share (\%) & 28.50 & 35.20 & +6.70 & +23.51\% \\
\bottomrule
\end{tabular}
\end{adjustbox}
\end{table}

\subsection{Sensitivity analysis}
\label{subsec:sensitivity}

Sensitivity analysis is conducted for two complementary objectives: (i) to identify which parameters most strongly influence the total-cost objective, and (ii) to identify which parameters most strongly affect feasibility under large deviations. The second objective is particularly important for multi-carrier energy hubs because infeasibility typically arises when deviations in renewable availability or demand exceed the combined flexibility of imports, dispatchable conversion, storage, and demand response.

A solver-based one-at-a-time (OAT) design was used. Each candidate parameter was perturbed over a predefined range while holding other inputs fixed, and the model was re-optimized for each perturbation. Sensitivity to the objective was quantified by the resulting change in total cost (and associated changes in emissions and carrier mix). Feasibility sensitivity was assessed by progressively widening perturbations until the model became infeasible, and recording which constraints were first violated.

The cost sensitivity results were summarized using normalized sensitivity indices and tornado plots. The parameters with the greatest impact on total cost were interpreted as the primary economic drivers of the optimal plan. Typical candidates include carbon price level, natural gas price, electricity import prices, electrolyzer CAPEX/O\&M, hydrogen storage CAPEX/O\&M, and demand levels. Here, demand response (DR) refers to a program where consumers adjust (shift) their electricity/heat use in response to incentives or prices, and the DR max shift ratios specify the maximum fraction of baseline load that can be moved up or down in any time period under that program.

Feasibility stress tests identify the parameters for which large deviations cannot be absorbed by the available flexibility. These parameters are interpreted as operationally critical uncertainties because moderate-to-large deviations may threaten feasibility (e.g., inability to satisfy electricity/heat/hydrogen-serving demands within import limits and storage bounds). In practice, feasibility-critical parameters frequently include renewable availability profiles, peak electricity and heat demands, and transport-related charging/refueling requirements, especially when they coincide with constrained imports or limited storage headroom.

We construct the uncertainty set in the robust optimization model of section \ref{subsec:robust} based on the feasibility sensitivity results. Specifically, parameters that (a) materially influence feasibility under stress, and (b) exhibit intrinsic short-term uncertainty are selected as uncertain inputs. This approach focuses robustness where it is operationally valuable and avoids unnecessary conservatism in protecting low-impact parameters.

\begin{figure}[ht]
\centering
\includegraphics[width=\linewidth]{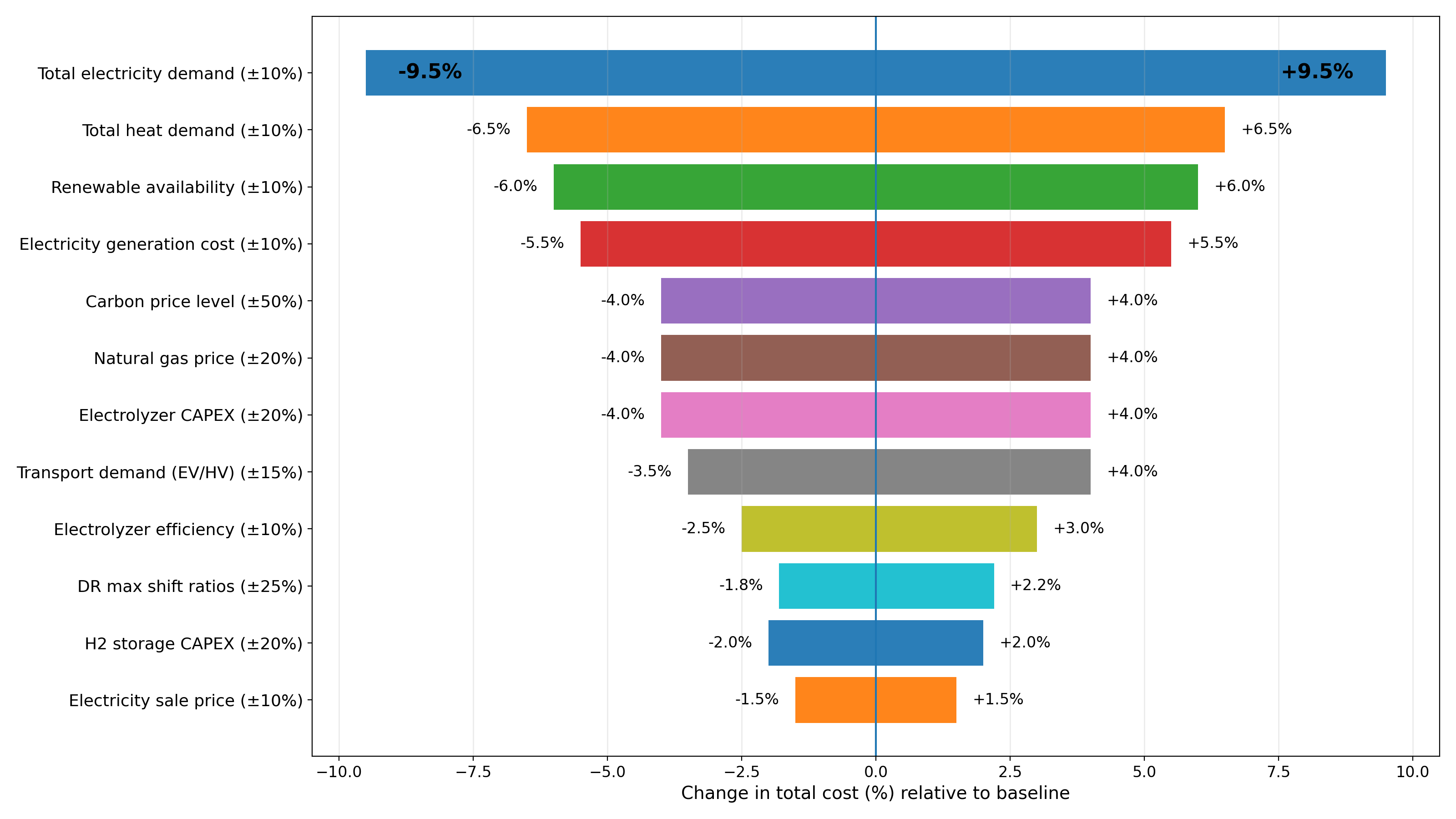}
\caption{Sensitivity tornado: one-at-a-time (OAT) perturbations showing the relative influence of candidate parameters on total cost. \\Bars indicate the range of cost changes under the low and high perturbation levels.}
\label{fig:tornado_cost}
\end{figure}

\subsection{Robust optimization model}
\label{subsec:robust}

\subsubsection{Method: budgeted robust optimization}
To manage operational decisions in uncertain environments while maintaining the tractability of MILP, this study uses the budgeted robust optimization framework from \cite{bertsimas2004price}. The main idea is to represent uncertain coefficients as bounded deviations around their expected (forecast) values. We protect key constraints against limited simultaneous worst-case deviations using a tunable conservatism parameter, known as the uncertainty budget.

Consider a general linear constraint $i$ that has some uncertain coefficients. Let $J_i$ represent the set of uncertain coefficients in constraint $i$. Each of these uncertain coefficients is assumed to lie within a symmetric interval around its expected value.
\begin{align}
\tilde{a}_{ij} = a_{ij} + \hat{a}_{ij} z_{ij}, \qquad z_{ij}\in[-1,1], \quad j\in J_i ,
\end{align}
where $a_{ij}$ is the expected coefficient, $\hat{a}_{ij}$ is the maximum deviation, and $z_{ij}$ is an unknown normalized deviation. If all $z_{ij}$ were allowed to take their worst-case values at the same time, the resulting solution would be overly conservative, which is equivalent to the classical ``full worst-case'' robustification. This approach addresses this condition by introducing a budget of uncertainty $\Gamma_i$ for each protected constraint $i$:
\begin{align}
\sum_{j\in J_i} |z_{ij}| \le \Gamma_i, \qquad 0 \le \Gamma_i \le |J_i|.
\end{align}

This budget specifies how many uncertain terms can be near their worst-case values simultaneously under one constraint. If $\Gamma_i=0$, then the robust model reduces to the expected model. If $\Gamma_i=|J_i|$, then the model will be nearly completely worst-case protected. For a linear constraint given as $\sum_j \tilde{a}_{ij} x_j \le b_i$, the robust counterpart is given as:
\begin{align}
\sum_{j} a_{ij} x_j \;+\; \max_{\{z_{ij}\}}\left\{\sum_{j\in J_i} \hat{a}_{ij} z_{ij} x_j:\; \sum_{j\in J_i} |z_{ij}| \le \Gamma_i,\; |z_{ij}|\le 1 \right\} \le b_i,
\end{align}
which ensures the feasibility of the solution in all the scenarios included in the budget. The robust constraint can be rewritten as linear constraints using auxiliary variables, thus allowing the linear robust counterpart to be incorporated into the MILP problem. The transformation introduces non-negative variables (often denoted as $p_i$ and $q_{ij}$), representing the budgeted deviation and the deviation in each element, and linear constraints to ensure the variables are within the limits set by the uncertain contributions. This approach maintains computational efficiency while allowing the control of the degree of conservatism through the budget parameter $\Gamma_i$.

In this study, the uncertainties are represented as symmetric $\pm 30\%$ deviation limits around the expected values, representing the short-term operating parameters that are immediately affected and have the most significant influence on the feasibility of the balance and coupling constraints. The uncertainties affect (i) the availability profiles of renewable energy sources, wind and solar power, and (ii) the time-varying demands, electricity and transport-related loads, including the electricity demand and transport-related loads, EV charging and hydrogen refueling, among others. The robustness budget parameter, $\Gamma$, is applied to control the degree of worst-case deviation throughout the 25-year planning horizon. The problem is only robustified, and robustness is applied only to the constraints involving uncertain terms, electricity balance, and hydrogen subsystem balance and supply constraints. The solution is thus guaranteed to be feasible under the worst-case scenario within the specified bounds, although at the cost of the associated robustness premium. Algorithm \ref{alg:bs_robust} summarizes the explained approach.

\begin{algorithm}[ht]
\caption{Budgeted robust optimization for the energy hub MILP}
\label{alg:bs_robust}
\DontPrintSemicolon

\KwIn{Nominal parameters $\bar{\theta}$; uncertainty bounds $\hat{\theta}$ (set as $\hat{\theta}=0.30\,\bar{\theta}$ for selected inputs); single robustness budget $\Gamma$; protected constraint set $\mathcal{I}^R$; MILP structure (variables, deterministic constraints, objective).
}
\KwOut{Robust solution $(x^{Rob},J^{Rob})$ and performance metrics including cost, emissions, fossil use, and H$_2$ production.}

\BlankLine
\textbf{Step 1: Identify uncertainty and protected constraints}\;
Select uncertain operational drivers $\Theta^{unc}$ (e.g., renewable availability and time-varying demands).\;
Identify the subset of constraints $\mathcal{I}^R$ in which $\Theta^{unc}$ appears (balance/coupling constraints).\;

\BlankLine
\textbf{Step 2: Define budgeted uncertainty set (Bertsimas--Sim)}\;
For each protected constraint $i \in \mathcal{I}^R$, define the uncertain coefficient index set $J_i$.\;
Represent each uncertain coefficient by $\tilde{a}_{ij} = a_{ij} + \hat{a}_{ij} z_{ij}$ with $0 \le z_{ij} \le 1$.\;
Impose the budget limit $\sum_{j\in J_i} z_{ij} \le \Gamma$ (single global budget).\;

\BlankLine
\textbf{Step 3: Build robust counterparts (linear reformulation)}\;
Initialize the robust MILP with the original objective and all deterministic constraints.\;
\ForEach{$i \in \mathcal{I}^R$}{
Replace constraint $i$ with its robust counterpart by introducing auxiliary variables $(p_i, q_{ij})$.\;
Add linear constraints that enforce $p_i + q_{ij} \ge \hat{a}_{ij}\,|x_j|$ for all $j\in J_i$ (using standard linearization with $y_j \ge |x_j|$).\;
Add the robustified constraint: $\sum_j a_{ij}x_j + \Gamma p_i + \sum_{j\in J_i} q_{ij} \le b_i$.\;
}

\BlankLine
\textbf{Step 4: Solve the robust MILP}\;
Solve the resulting robust MILP using CPLEX in GAMS.\;
Extract robust decision variables and compute scenario metrics.\;

\end{algorithm}

\subsubsection{Robust results under carbon tax and net zero scenarios}

Robust outcomes at the Scenario-level are summarized in Table~\ref{tab:scenario_metrics}. These results are calibrated using a present-year (current) baseline dataset for operational inputs (e.g., demand levels, renewable availability profiles, and price/cost parameters), rather than externally projected future time-series. The long-term horizon to 2050 is therefore interpreted as a planning evaluation under current baseline conditions and policy assumptions, with robustness capturing the impact of short-term uncertainty around those baseline inputs.

\begin{figure}[H]
\centering
\includegraphics[width=.8\linewidth]{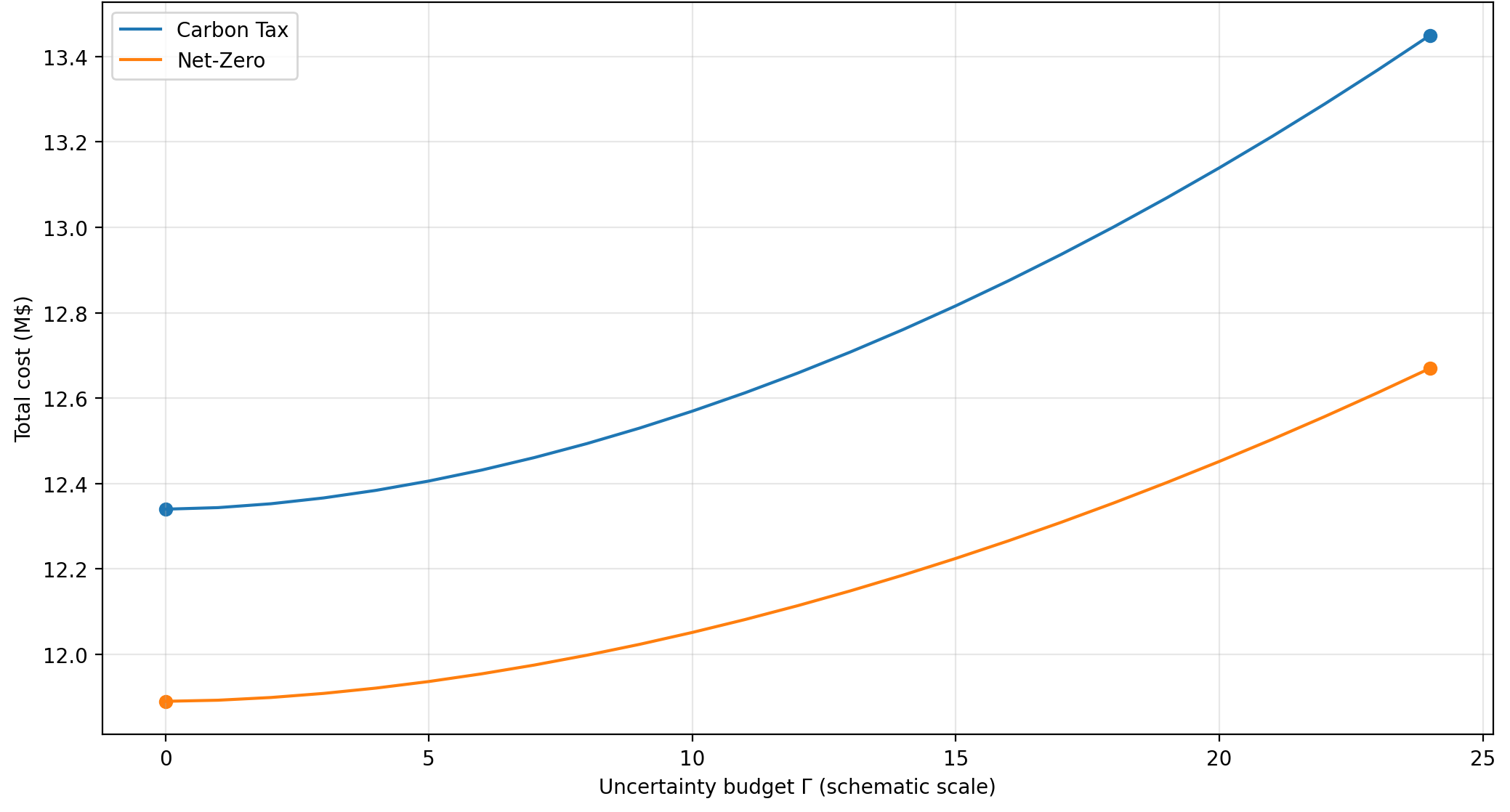}
\caption{Robustness premium as the uncertainty budget $\Gamma$ increases (2025 baseline). \\Larger $\Gamma$ implies protection against more simultaneous adverse deviations, increasing the objective value.}
\label{fig:cost_vs_gamma}
\end{figure}

In both policy scenarios, the addition of robustness results in increased costs compared to the deterministic case. In the carbon tax scenario, the costs increase from 12.34 to 13.45 million dollars, a 9.0\% growth. In the net zero scenario, costs increase from 11.89 to 12.67 million dollars, a growth of 6.56\%. Figure \ref{fig:cost_vs_gamma} below demonstrates this process, where an increased level of uncertainty, $\Gamma$, is accommodated while maintaining feasibility through a higher number of concurrent adverse scenarios. This results in additional operational costs, increasing the objective value.

The addition of robustness also results in changes in the level of energy carrier usage. In the carbon tax scenario, CO$_2$ emissions increase from 0.45 to 0.52 million tons, a 15.6\% growth. The use of fossil fuels increases from 25.6 to 28.3 GWh/year, a growth of 10.5\%. Hydrogen production decreases from 15.2 to 14.8 GWh/year, a 2.6\% decrease. Electrolyzer utilization decreases slightly to 60.8\%, down by 1.5\%, while the renewable share decreases slightly to 27.9\%, down by 0.6\%. The overarching conclusion here is that robustness maintains dispatchability flexibility in the presence of uncertainty, which in this case is highly dependent on fossil fuels. This is exemplified by electrolyzer scheduling during stress scenarios.

A similar trend can be seen in the net zero policy. Total costs increase from 11.89 to 12.67 million dollars. Emissions increase from 0.39 to 0.44 million tons of CO$_2$, up by 12.8\%. Hydrogen production decreases slightly from 18.5 to 17.9 GWh/year, down by 3.2\%. Fossil fuel use increases from 19.18 to 23.06 GWh/year, up by 20.2\%. Meanwhile, the renewable share and the utilization of electrolyzers decrease slightly by 1\% and 1.5\% to 34.6\% and 73.9\%, respectively. In conclusion, robust solutions give up a certain level of benefits in terms of decarbonization in favor of feasibility under uncertainty within the specified range of $\pm30\%$ and the chosen budget $\Gamma$.

\begin{table}[ht]
\centering
\caption{Scenario-level performance metrics for deterministic (Det) and robust (Rob) solutions under carbon tax and net zero policy cases.}
\label{tab:scenario_metrics}
\begin{adjustbox}{width=\textwidth,center}
\begin{tabular}{lccccc}
\toprule
\textbf{Metric} & \textbf{Unit} & \textbf{Carbon Tax (Det)} & \textbf{Carbon Tax (Rob)} & \textbf{Net Zero (Det)} & \textbf{Net Zero (Rob)} \\
\midrule
Total Cost & M\$ & 12.34 & 13.45 & 11.89 & 12.67 \\
Total Emissions & MTon CO$_2$ & 0.45 & 0.52 & 0.39 & 0.44 \\
H$_2$ Production & GWh/year & 15.2 & 14.8 & 18.5 & 17.9 \\
Fossil Fuel Use & GWh/year & 25.6 & 28.3 & 19.18 & 23.06 \\
Renewable Share & \% & 28.5 & 27.9 & 35.2 & 34.6 \\
Electrolyzer Utilization & \% & 62.3 & 60.8 & 75.4 & 73.9 \\
\bottomrule
\end{tabular}
\end{adjustbox}
\end{table}

\subsubsection{Deterministic versus robust comparison and managerial interpretation}

The comparison between deterministic and robust results provides quantification of the robustness premium and helps identify the mechanisms through which the model hedges uncertainty. First, the total cost increases by 6.6\% to 9.0\% in all policy cases (Table~\ref{tab:scenario_metrics}), as the buffers are increased to ensure the feasibility of the solution under uncertainty. Second, the model's robustness solution shows increased emissions and fossil fuel consumption in all policy cases, accompanied by slightly reduced production of hydrogen, renewable energy, and the use of the electrolyzer. Overall, the results suggest that when faced with uncertainty, the reliance on more balancing through dispatchable energy sources increases, especially where the uncertainty is related to balance and coupling constraints.

From the perspective of the planner, the results indicate the areas where investment and operational policy can mitigate the emissions implications of the robustness solution while maintaining the solution's reliability. The implications of the robust solution are that the existing flexibility capacity is not sufficient to hedge the uncertainty, and the increased fossil fuel balancing is necessary. By adding more flexibility options that are non-emitting (such as additional storage capacity, increased demand response, renewable energy firming, and/or dispatchable non-emitting energy), the model's ``hedging'' can be reduced through fossil and emissions-intensive balancing. From the managerial perspective, the implications are that the concept of robustness is not merely an ``extra cost,'' but rather an opportunity to identify the areas where the system lacks flexible capacity and where investment policy can reduce the robustness premium and the emissions implications of the robust solution.

To further illustrate the process of propagation of uncertainty in the robustness premium, two parameters will be discussed: total market demand, and hydrogen energy price. These parameters were chosen because they represent two different mechanisms. Uncertainty in total market demand directly impacts balance constraints related to feasibility-critical variables such as electricity, heat, and hydrogen service. It is considered a main driver for buffers in the robust solution. In contrast, uncertainty in the price of hydrogen primarily affects the economic attractiveness of hydrogen pathways (electrolyzer operation and storage use), changing the cost-optimal mix without necessarily tightening the feasibility to the same extent.

\begin{figure}[H]
\centering
\begin{subfigure}[t]{0.49\linewidth}
    \centering
    \includegraphics[width=\linewidth]{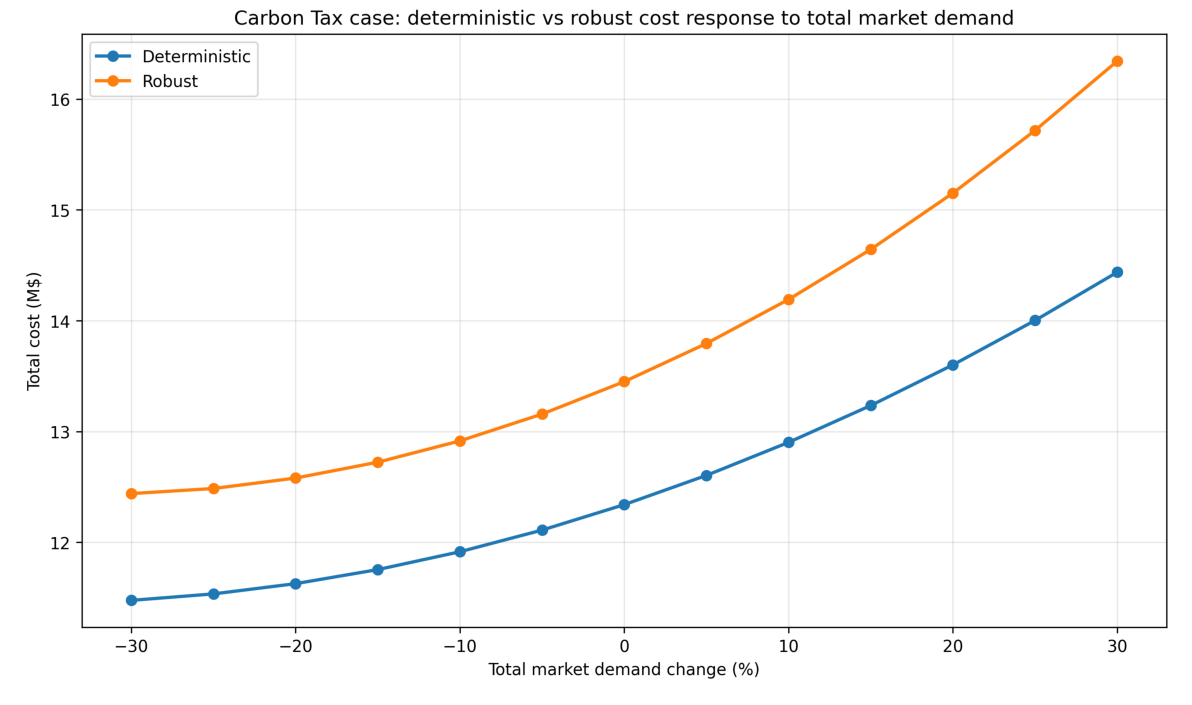}
    \caption{Total market demand (Det vs Rob).}
    \label{fig:ct_detrob_demand}
\end{subfigure}
\hfill
\begin{subfigure}[t]{0.49\linewidth}
    \centering
    \includegraphics[width=\linewidth]{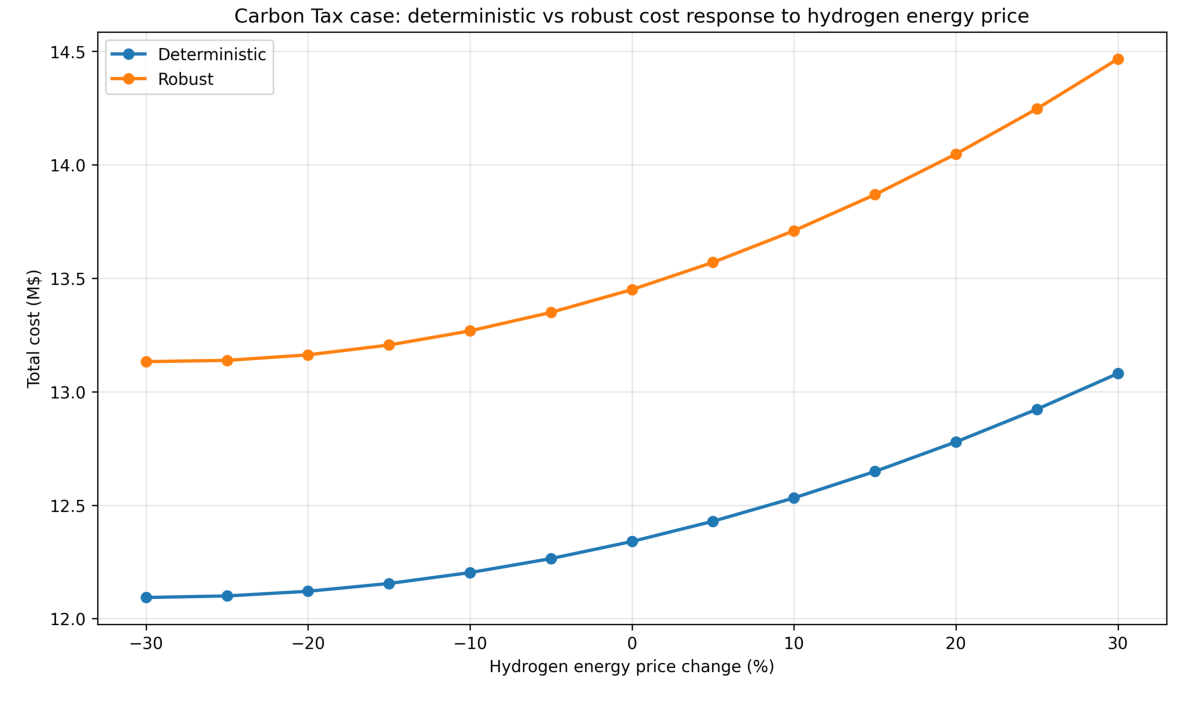}
    \caption{Hydrogen energy price (Det vs Rob).}
    \label{fig:ct_detrob_h2price}
\end{subfigure}
\caption{Carbon tax case (2025 baseline): deterministic vs robust total cost response under one-at-a-time perturbations in (a) total market demand and (b) hydrogen energy price.}
\label{fig:ct_detrob_two_params}
\end{figure}

\begin{figure}[H]
\centering
\begin{subfigure}[t]{0.49\linewidth}
    \centering
    \includegraphics[width=\linewidth]{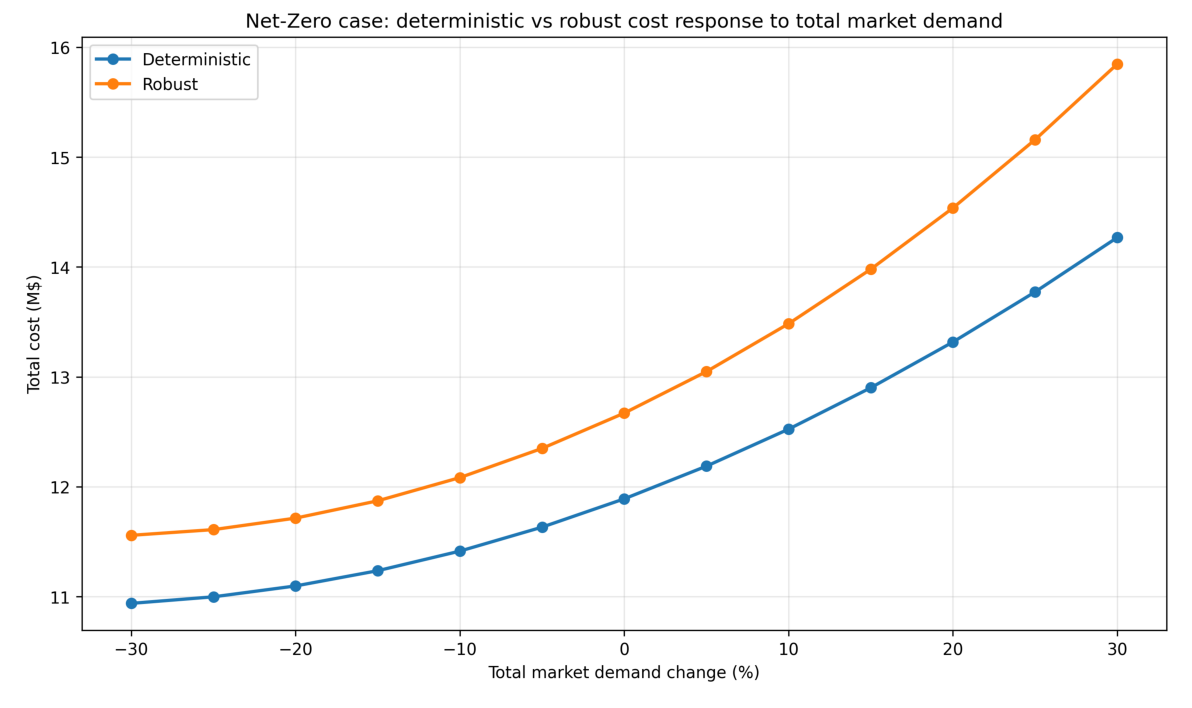}
    \caption{Total market demand (Det vs Rob).}
    \label{fig:nz_detrob_demand}
\end{subfigure}
\hfill
\begin{subfigure}[t]{0.49\linewidth}
    \centering
    \includegraphics[width=\linewidth]{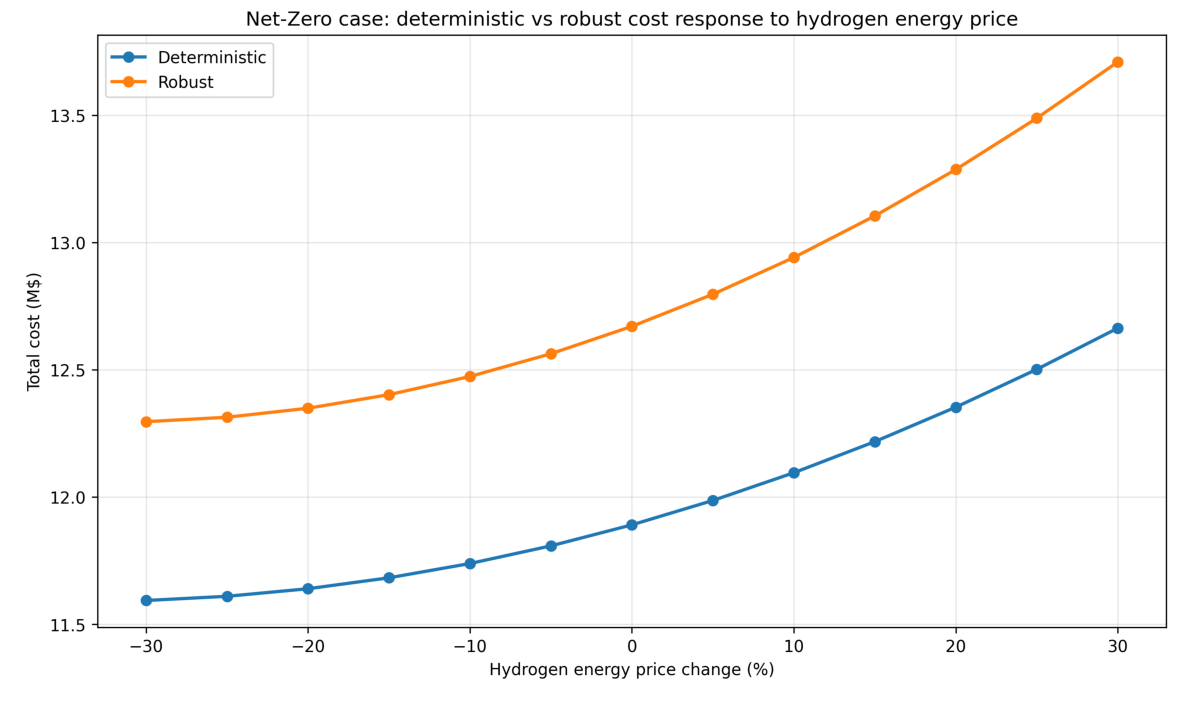}
    \caption{Hydrogen energy price (Det vs Rob).}
    \label{fig:nz_detrob_h2price}
\end{subfigure}
\caption{Net zero case (2025 baseline): deterministic vs robust total cost response under one-at-a-time perturbations in (a) total market demand and (b) hydrogen energy price.}
\label{fig:nz_detrob_two_params}
\end{figure}

Figures \ref{fig:ct_detrob_two_params} and \ref{fig:nz_detrob_two_params} illustrate the variation in total costs when the total market demand changes. In all cases, it can be seen that the robust cost curves remain above the deterministic curves over the range of perturbations. These plots illustrate the robustness premium to guaranty feasibility under uncertain conditions. Furthermore, increasing total market demands inflates the difference between robust and deterministic solutions. This underlines the impact of total market demand on balance constraints, where growing demand directly stresses flexibility buffers.

Although often less pronounced than the impact of total market demand, the price of hydrogen energy influences both deterministic and robust costs. This implies that this uncertainty impacts the economic attractiveness of different service paths, such as electrolyzer operations and storage utilization, while total market demand directly jolts the flexibility buffers required to maintain feasibility.

\begin{figure}[ht]
\centering
\includegraphics[width=\linewidth]{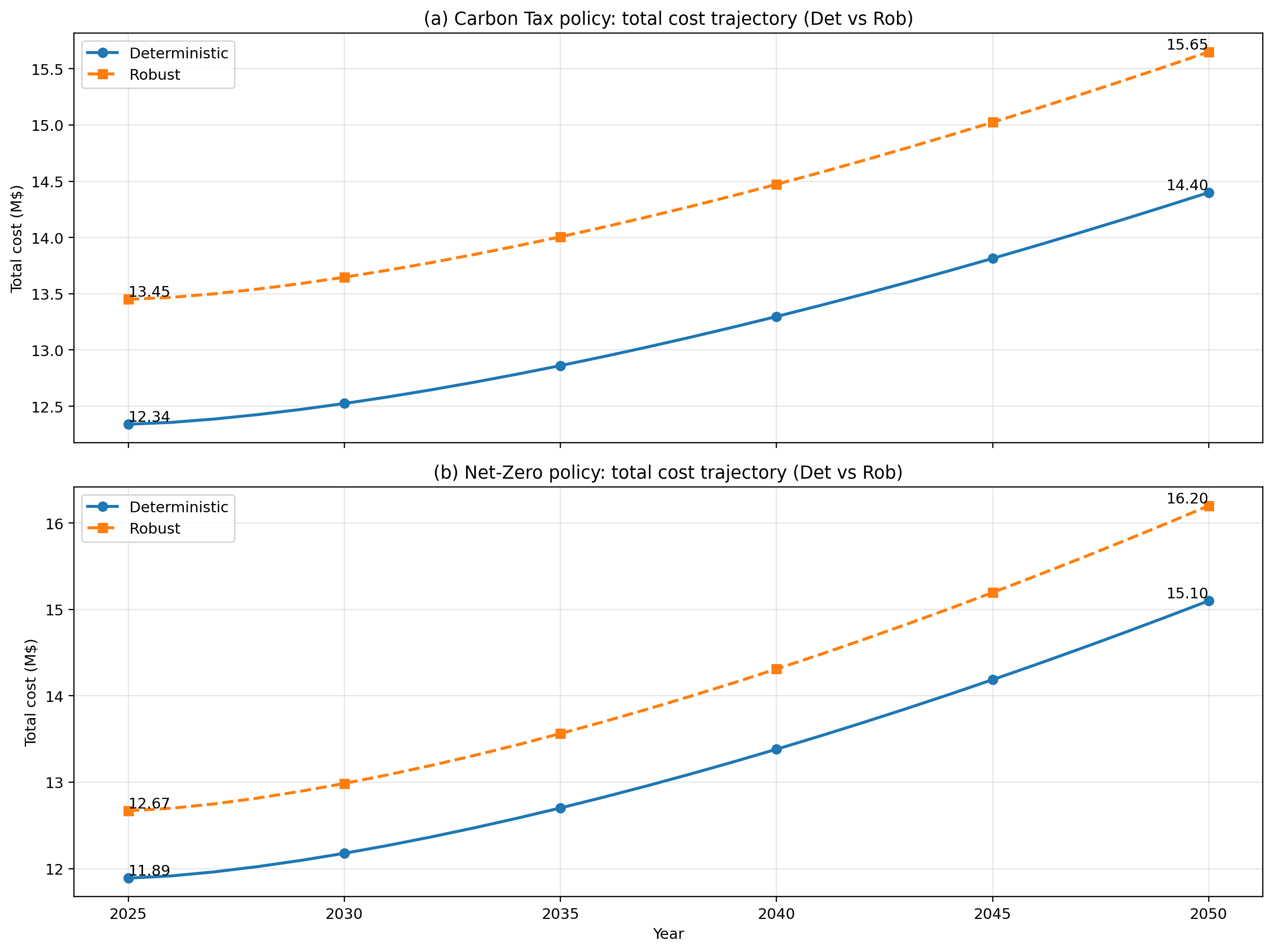}
\caption{Total cost trajectories (2025--2050): carbon tax policy (top) and net zero policy (bottom), comparing deterministic and robust solutions.}
\label{fig:cost_trajectories_policy_det_rob}
\end{figure}

Figure \ref{fig:cost_trajectories_policy_det_rob} compares the total system cost trajectories under the two carbon abatement schemes when solved deterministically and using budgeted robust optimization. Note that in the two plots, the robust solution is always more costly than the deterministic solution, reflecting the cost of robustness, i.e., the ``robustness premium'' that must be paid to ensure the solution's feasibility in the face of bounded adverse changes in the parameters selected as uncertain. The gap between the two solution methods grows over time, indicating that as the system grows and becomes more constrained, the cost of additional buffers becomes more pronounced and more expensive.

\section{Conclusions}\label{sec:CC}

This research proposes a robust carbon-conscious optimization approach for incorporating hydrogen into a grid-connected multi-energy hub that offers electric, thermal, and transportation energy services. The basic model is an MILP that includes balance and coupling relations for multiple energy carriers, conversion, storage, and demand response \citep{Siroos2025}. First, the basic model is expanded to a deterministic model, which is used to determine the influential uncertainty parameters. Afterward, a budgeted robust optimization approach is introduced to address the operational and feasibility-related uncertainties identified, with $\pm 30\%$ perturbation for the uncertain variables and a single uncertainty budget $\Gamma$. The research assumes a 25-year planning horizon (2025--2050) and contemplates two different approaches to decarbonization, i.e., carbon tax and net zero policies. The research selected the province of Ontario, Canada as a real-world case study and designated a cold and a warm day to represent the hub operation.

The deterministic results show that hydrogen is a viable alternative energy carrier and a feasible buffering service for a power, heat, and transportation hub. The results of the long-term planning scenario show that the model predicts a significant expansion of the hydrogen infrastructure through 2050, with the electrolyzer capacity increasing to 3,800 MW by 2050 compared to 300 MW in 2025, and the hydrogen storage capacity increasing to 37,000 MWh compared to 2,000 MWh, while also showing a significant increase in hydrogen production. The results for the policy scenario show that the net zero constraint results in a stronger reduction of emissions and fossil fuels compared to the carbon tax scenario while also increasing hydrogen production and renewable use, as technology substitution is stronger for an emissions trajectory constraint.

Robust optimization outcomes provide quantified measures of the cost implications of feasibility protection under uncertainty. In terms of overall cost, the robust approach prioritizes over its deterministic counterpart in both carbon tax and net zero cases, providing a robustness premium between 6.6\% and 9.0\%. In terms of overall outcomes, robust optimization leads to increased emissions and fossil fuel consumption, as well as a small decrease in hydrogen production, electrolyzer utilization, and renewable penetration. This suggests that, under bounded worst-case uncertainty, the system retains dispatchability to ensure feasibility, and limited low-carbon flexibility allows for robustness to be partly mitigated by fossil-based flexibility. From a planning and managerial point of view, these results demonstrate opportunities to reduce emissions without impacting robustness, which may include increased H$_2$ storage capacity, enhanced demand response programs, and increased renewable firming, among others.

Several avenues exist for further improvements in future works. First, the operating model currently utilizes representative cold and warm day operating horizons. Expanding this to a more detailed seasonal or even a meticulous chronological operating sequence would further enhance the model rigor. Additionally, this study models electric vehicle and fuel cell electric vehicle demands at an aggregate level. Future research could expand this approach to detailed mobility and fueling constraints. Furthermore, we suggest that future studies utilize a hybrid stochastic-robust approach that combines tractable feasibility protection with distributional information. Finally, the inclusion of more detailed network constraints, equipment operating details, and technology learning would further enhance the scalability of the developed model.

\bibliographystyle{elsarticle-harv}
\bibliography{References}

\clearpage

\appendix
\section{Notations}
\label{Appendix:Notations}

\setcounter{table}{0}

\small
\renewcommand{\arraystretch}{0.7}
\begin{longtable}{@{}p{0.24\textwidth}p{0.71\textwidth}@{}}
\label{tab:notation_ontario}\\
\toprule
\textbf{Symbol} & \textbf{Description} \\
\midrule
\endfirsthead

%\multicolumn{2}{c}{\tablename\ \thetable\ -- continued from previous page} \\
\toprule
\textbf{Symbol} & \textbf{Description} \\
\midrule
\endhead

\midrule
\multicolumn{2}{r}{Continued on next page} \\
\endfoot

\bottomrule
\endlastfoot

\multicolumn{2}{@{}l}{\textbf{Sets and indices}}\\
$\mathcal{Y}$ & Set of years in the planning horizon, indexed by $y$.\\
$\mathcal{T}$ & Set of intra-year operating periods, indexed by $t$.\\
$\mathcal{I}^{R}$ & Set of constraints protected in the robust counterpart, indexed by $i$.\\
$J_i$ & Index set of uncertain coefficients appearing in protected constraint $i$.\\
$\Theta^{\mathrm{unc}}$ & Set of uncertain operational drivers considered in the robust model.\\
$\chi \in \{el,h\}$ & Energy-demand category used in demand-response notation (electricity or heat).\\[5mm]

\multicolumn{2}{@{}l}{\textbf{Objective and aggregate cost terms}}\\
$J_{\mathrm{cost}}$ & Total operating cost objective over the planning horizon.\\
$J_{\mathrm{emis}}$ & Total direct-emissions objective over the planning horizon.\\
$C_{\mathrm{cyc}}(y)$ & Storage-cycling cost in year $y$.\\
$C_{\mathrm{OM}}(y)$ & Operation and maintenance cost in year $y$.\\
$C_{\mathrm{EM}}(y)$ & Emissions monetization term in year $y$.\\
$C_{\mathrm{DR}}(y)$ & Demand-response cost in year $y$.\\
$C_{\mathrm{EDR}}(y),\, C_{\mathrm{HDR}}(y)$ & Electricity- and heat-demand-response cost components in year $y$.\\[5mm]

\multicolumn{2}{@{}l}{\textbf{Parameters}}\\
$\pi^{\mathrm{buy}}_{y,t},\, \pi^{\mathrm{sell}}_{y,t}$ & Electricity import and export prices in year $y$ and period $t$.\\
$P_f$ & Emission factor or emissions-cost coefficient applied to fuel use.\\
$\tau_{\mathrm{CO2}}$ & Carbon-tax parameter used in the carbon-pricing policy case.\\
$\mathrm{LHV}$ & Lower heating value of the relevant fuel.\\
$L^{el}_{y,t},\, L^{h}_{y,t}$ & Baseline electricity and heat demand in year $y$ and period $t$.\\
$P^{HV}_{y,t}$ & Hydrogen-vehicle demand in year $y$ and period $t$.\\
$MR^{\uparrow}_{\chi},\, MR^{\downarrow}_{\chi}$ & Maximum upward and downward demand-response shift ratios for demand type $\chi$.\\
$\eta^{\mathrm{B}}$ & Gas-boiler efficiency.\\
$\eta^{\mathrm{CHP}}_e,\, \eta^{\mathrm{CHP}}_h$ & Electrical and thermal efficiencies of the CHP unit.\\
$\eta^{\mathrm{ely}}$ & Electrolyzer efficiency.\\
$\eta^{\mathrm{ch}}_{\mathrm{H}},\, \eta^{\mathrm{dch}}_{\mathrm{H}}$ & Thermal-storage charging and discharging efficiencies.\\
$\eta^{\mathrm{ch}}_{\mathrm{H2}},\, \eta^{\mathrm{dch}}_{\mathrm{H2}}$ & Hydrogen-storage charging and discharging efficiencies.\\
$H^{\min},\, H^{\max}$ & Minimum and maximum thermal-storage state.\\
$H2^{\min},\, H2^{\max}$ & Minimum and maximum hydrogen-storage state.\\
$P^{\max}_{\mathrm{ch,H}},\, P^{\max}_{\mathrm{dch,H}}$ & Maximum thermal-storage charging and discharging powers.\\
$P^{\max}_{\mathrm{ch,H2}},\, P^{\max}_{\mathrm{dch,H2}}$ & Maximum hydrogen-storage charging and discharging powers.\\
$P^{\max}_{\mathrm{ely}}$ & Nameplate capacity of the electrolyzer.\\
$P^{\max}_r$ & Ramp-rate limit for a ramp-constrained generator.\\
$\bar{\theta},\, \hat{\theta}$ & Nominal values and maximum deviations of selected uncertain inputs in the robust model.\\[5mm]

\multicolumn{2}{@{}l}{\textbf{Continuous operating variables and state variables}}\\
$P^{\mathrm{buy}}_{y,t},\, P^{\mathrm{sell}}_{y,t}$ & Imported and exported electrical power of the hub in year $y$ and period $t$.\\
$P^{\mathrm{PV}}_{y,t},\, P^{\mathrm{WT}}_{y,t}$ & Electrical power produced by the photovoltaic and wind units.\\
$P^{\mathrm{ECHP}}_{y,t}$ & Electrical power generated by the CHP unit.\\
$P^{\mathrm{EL}}_{y,t}$ & Native electrical load served by the hub.\\
$P^{\mathrm{EV}}_{y,t}$ & Electric-vehicle load served by the hub.\\
$P^{\mathrm{ice}}_{y,t}$ & Electrical power allocated to the cooling-storage/ice subsystem.\\
$P^{\mathrm{ec}}_{y,t}$ & Electrical power consumed by the electric chiller.\\
$P^{\mathrm{ely}}_{y,t}$ & Electrical power consumed by the electrolyzer.\\
$P^{\mathrm{HB}}_{y,t}$ & Heat produced by the gas boiler.\\
$P^{\mathrm{HCHP}}_{y,t}$ & Heat produced by the CHP unit.\\
$P^{\mathrm{Hac}}_{y,t}$ & Heat supplied to the absorption chiller.\\
$G^{\mathrm{B}}_{y,t}$ & Fuel consumed by the gas boiler.\\
$G^{\mathrm{CHP}}_{y,t}$ & Fuel consumed by the CHP unit.\\
$G^{\mathrm{gas}}_{y,t},\, G^{\mathrm{bio}}_{y,t}$ & Gas and biofuel use contributing to direct emissions.\\
$H_{y,t}$ & Thermal energy stored in the heat-storage  at the end of period $t$.\\
$P^{\mathrm{ch,H}}_{y,t},\, P^{\mathrm{dch,H}}_{y,t}$ & Thermal-storage charging and discharging powers.\\
$H2_{y,t}$ & Hydrogen inventory at the end of period $t$.\\
$P^{\mathrm{ch,H2}}_{y,t},\, P^{\mathrm{dch,H2}}_{y,t}$ & Hydrogen-storage charging and discharging powers.\\
$P^{\mathrm{H2,use}}_{y,t}$ & Hydrogen supplied from storage to satisfy hydrogen-service demand.\\
$P^{\uparrow}_{\chi,y,t},\, P^{\downarrow}_{\chi,y,t}$ & Upward and downward demand-response adjustments for demand type $\chi$.\\
$P_{y,t}$ & Generic output of a ramp-limited generation unit in period $t$.\\[5mm]

\multicolumn{2}{@{}l}{\textbf{Binary variables}}\\
$K^{\mathrm{ch}}_{y,t},\, K^{\mathrm{dch}}_{y,t}$ & Binary variables for thermal-storage charging and discharging modes.\\
$K^{\mathrm{ch,H2}}_{y,t},\, K^{\mathrm{dch,H2}}_{y,t}$ & Binary variables for hydrogen-storage charging and discharging modes.\\
$I^{\uparrow}_{\chi,y,t},\, I^{\downarrow}_{\chi,y,t}$ & Binary variables indicating upward and downward demand-response activation for demand type $\chi$.\\[2mm]

\multicolumn{2}{@{}l}{\textbf{Robust-optimization auxiliary notation}}\\
$\tilde{a}_{ij}$ & Uncertain coefficient in protected constraint $i$.\\
$a_{ij}$ & Nominal value of uncertain coefficient $\tilde{a}_{ij}$.\\
$\hat{a}_{ij}$ & Maximum deviation of uncertain coefficient $\tilde{a}_{ij}$.\\
$z_{ij}$ & Normalized deviation variable associated with coefficient $\tilde{a}_{ij}$.\\
$\Gamma_i$ & Uncertainty budget for protected constraint $i$.\\
$\Gamma$ & Global robustness-budget parameter used in the budgeted-uncertainty model.\\
$b_i$ & Right-hand-side value of protected constraint $i$.\\
$x_j$ & Generic decision variable in the robust counterpart.\\
$y_j$ & Auxiliary nonnegative variable used to linearize $|x_j|$ when needed.\\
$p_i,\, q_{ij}$ & Auxiliary variables used in the linear reformulation of the Bertsimas--Sim robust counterpart.\\

\end{longtable}
\normalsize

\end{document}